\documentstyle[leqno,12pt]{article}
\begin{document}
\baselineskip=0.7cm
\textheight=9.0in

\begin{titlepage}
\begin{center}
     {\bf \Large Convergence Rates of Finite Difference Stochastic Approximation Algorithms \footnote{\em This work is supported in part by the U.S. Army Research Office under agreement W911NF-04-D-0003.}}
\end{center}

\begin{center}
        Liyi Dai \\
%        Computing Sciences Division\\
        Army Research Office \\
        Research Triangle Park, NC 27703\\
        liyi.dai.civ@mail.mil \\[4mm]
\end{center}

\vspace{0.5cm}

\centerline{\bf Abstract}
%\vspace{0.7cm}
\begin{center}
\parbox{5.5in}{\baselineskip=0.7cm
Recently there has been renewed interests in derivative free approaches to stochastic 
optimization. In this paper, we examine the rates of convergence for the 
Kiefer-Wolfowitz algorithm and the mirror descent algorithm, under various updating 
schemes using finite differences as gradient approximations. It is shown that the 
convergence of these algorithms can be accelerated by controlling the implementation 
of the finite differences. Particularly, it is shown that the rate can be increased to
 $n^{-2/5}$ in general and to $n^{-1/2}$ in Monte Carlo optimization for a broad 
 class of problems, in the iteration number $n$. }
\end{center}

\vspace{0.5cm}

\vspace{0.5cm}
%{\bf AMS subject classifications.} 90C15, 90C25
\vspace{0.5cm}

\end{titlepage}

\topmargin=-0.5in

{\bf 1. Introduction.}\hspace{2mm}
Let $R$ denote the set of real numbers. 
%For a real-valued differentiable
%function $f(x,y)$, we will use $f'_{x}(x,y)$ to denote its partial 
%derivative with respect to $x$. The subscript is omitted for functions 
%of one variable. 
Consider a real-valued function
$J(\theta)$ of the form $J(\theta)=E_{X}[L(X(\theta))]$ where $\theta$ is
a parameter, or a vector of parameters, $L(X)$ is a real-valued function,
and $X(\theta)$ is a random variable that depends on $\theta$. For 
simplicity, throughout this paper we assume that $\theta$ is a scalar
and $\theta \in \Theta \subset R$, and $X(\theta)$ is of the form 
$X(\theta)= X(\theta,\xi)$, where $\xi$ is a random variable independent 
of $\theta$. In such a formulation, $X(\theta)$ is parameterized on an 
underlying probability space that is independent of $\theta$. For any two 
random variables $\eta$ and $\xi$, there exists a Borel function $\phi$ 
such that $\eta=\phi(\xi)$ [Shiryayev (1984), p.172]. Such a representation for 
$X(\theta,\xi)$ is always possible. Therefore,
$J(\theta)$ can be written as $J(\theta) = 
E_{\xi}[L(X(\theta,\xi))]$. We are particularly interested in finding 
an optimal parameter $\theta^{*}\in \Theta$ to optimize, say minimize, 
$J(\theta)$. This is a challenging problem since the analytical form 
of $J(\theta)$ is usually unavailable for most problems of interest.
What is obtainable are the sample
measurements of the random value of $L(X(\theta,\xi))$. We have to use the 
information on $L(X(\theta,\xi))$ to find $\theta^{*}$. Such stochastic optimization 
problems can be found in many applications. The main
approach to finding the optimal solution is to successively approximate
$\theta^{*}$ via algorithms of {\em stochastic approximation}. This is a 
classical and standard approach that has been adopted in practice for 
decades. The {\em Robbins-Monro} (RM) algorithm, the 
{\em Kiefer-Wolfowitz} (KW) algorithm, and the relatively recent  
mirror descent (MD) algorithm  are the most popular algorithms
of this class.

The RM algorithm, introduced by Robbins and Monro (1951), finds $\theta^{*}$
in the following way. Let $\theta_{0}$ be selected and 
$\{a_{n}\}$ a sequence of positive numbers. For each integer 
$n\geq 0$, let
\begin{equation}
       \theta_{n+1}=\theta_{n}-a_{n}g_{n}
\label{rm}
\end{equation}
where $g_{n}$ is an unbiased estimate of the derivative $J'(\theta)$
of $J(\theta)$ with respect to $\theta$. Assume that $J'(\theta)$ exists 
on $\Theta$ and the variance of $g_{n}$ is uniformly bounded for all $n$. 
Assume $(\theta-\theta^{*})J'(\theta)>0$ for all $\theta\not= \theta^{*}$, 
\[  \sum_{n}a_{n}=\infty, \;\; \sum_{n}a_{n}^{2}<\infty,     \]
and that several other technical conditions are satisfied. Then $\{\theta_{n}\}$ 
converges to $\theta^{*}$ with probability one. The convergence rate 
(in terms of root mean square error) is $n^{-1/2}$. Note that this
is the best possible rate of convergence for algorithms of the form (1) for
stochastic optimization [see, e.g., Fabian (1971)].

The KW algorithm, introduced by Kiefer and Wolfowitz (1952), is a modification 
of the RM algorithm by approximating the gradient using a finite difference and  
finds $\theta^{*}$ recursively by
\begin{equation}
    \theta_{n+1} = \theta_{n}-a_{n}h_{n},
\label{kw}
\end{equation}
where
\begin{equation}
    h_{n}
         =\frac{L(X(\theta_{n}+\delta_{n},\xi_{1,n}))
          -L(X(\theta_{n}-\delta_{n},\xi_{2,n}))}{2\delta_{n}},  
\label{kw_h}
\end{equation}
$\{\delta_{n}\}$ is a sequence of positive numbers,
$L(X(\theta_{n}+\delta_{n},\xi_{1,n}))$ and $L(X(\theta_{n}-\delta_{n},
\xi_{2,n}))$ are two measurements of $L(X(\theta,\xi))$ at $\theta_{n}
+\delta_{n}$ and $\theta_{n}-\delta_{n}$, and $\xi_{1,n}, \xi_{2,n}$ are
corresponding samples of $\xi$. Kiefer-Wolfowitz (1952)
proved that if $J(\theta)$ is decreasing for $\theta<\theta^{*}$ and
increasing for $\theta>\theta^{*}$, and if
\[  \delta_{n}\rightarrow 0,\;\; \sum_{n}a_{n}=\infty, \;\;
      \sum_{n}a_{n}\delta_{n}<\infty, \;\; 
              \sum_{n}a_{n}^{2}/\delta_{n}^{2}<\infty, \]
the sequence $\{\theta_{n}\}$ converges to $\theta^{*}$ with 
probability one under some additional minor conditions. If all entries in
$\{\xi_{i,n}\}$ are mutually independent, the best possible convergence
rate for the KW algorithm (2) is $n^{-1/3}$ which is achieved by choosing 
$a_{n}=an^{-1}, \delta_{n}=d n^{-1/6}$ with $a, d > 0$
constants [e.g. Burkholder (1956); Fabian (1971); Sacks (1958)]. The rate 
$n^{-1/3}$ is regarded not satisfactory compared to the best 
possible rate $n^{-1/2}$ for the RM algorithm. 

The MD algorithm, introduced by Nemirovski and Yudin (1983), improves 
the robustness of gradient based optimization algorithms. At iteration $n\geq 0$, 
$\theta_{n+1}$ is updated via solving
\begin{equation}
    \theta_{n+1} = \textrm{argmin}_{\theta\in\Theta}\left\{ <h_n,\theta>+\frac{1}{a_n}D_{\psi}(\theta,\theta_n)\right\},
\label{md}
\end{equation}
where $h_n$ is an estimate of the derivative $J'(\theta)$, 
%\[  h_n = g_n\;\; \textrm{  or  } \;\;  h_n
%         =\frac{L(X(\theta_{n}+\delta_{n},\xi_{1,n}))
%          -L(X(\theta_{n},\xi_{2,n}))}{\delta_{n}}, \]  
$D: \Theta\times\Theta\rightarrow R^+$ 
is a Bregman distance defined as
\begin{equation}
    D(\theta,\tau) := \psi(\theta)-\psi(\tau)-<\psi'(\tau),\theta-\tau> \geq \kappa ||\theta-\tau||^2,
 \label{md_bd}
\end{equation}
where $\psi(.)$ is a distance generating function and $\kappa>0$ is a constant. 
In (\ref{md_bd}), $||.||$ is a general norm on $R^m$ (and on $R$ in this paper).
It has been established by Nemirovski et al. (2009) and Duchi et al. (2012,2013)
 that if $J(\theta)$ is convex, Lipschitz continuous and
\[ a_n \rightarrow 0, \;\; \sum_{n}a_{n}=\infty, \;\; \]
\[ \hat{\theta}_n = \frac{1}{n}\sum_{i=1}^n \theta_i \textrm{  or   } \hat{\theta}_n 
=\sum_{i=1}^n \nu_i\theta_i,\;\; \nu_i = \frac{a_i}{\sum_{j=1}^i a_j}, i =1 , 2, ... \]
then $J(\hat{\theta}_n)$ converges to the minimum of $J(\theta)$ and the rate of convergence is 
$n^{-1/2}$ under mild technical conditions that will be specified in Section 3.

%As a result of decades-long effort made by numerous researchers, the theory on 
The convergence of these algorithm is 
fairly understood [Burkholder (1956); Fabian (1971); Kushner 
and Clark (1978); Chung (1954); Dupa$\check{c}$ (1957); Dvoretzky (1956);
Sacks (1958); Wasan (1969), Nemirovski et al. (2009), Duchi et al. (2012, 2013)]. 
The conditions for the convergence of these 
algorithms can be made substantially weaker than those we have previously 
mentioned [e.g. Kushner and Clark (1978); Wasan (1969)]. 
The convergence rate for the RM algorithm is much faster than that for 
the KW algorithm. This is not surprising if we note that the KW 
algorithm uses the finite difference $h_{n}$ as an approximation 
to the derivative $J'(\theta)$, while the RM algorithm uses an unbiased 
estimate of $J'(\theta)$. Therefore, the faster rate is achieved 
at the cost of obtaining an unbiased estimate of the derivative
$J'(\theta)$ that is often challenging in practice.
%Although obtaining an estimate of the derivative is not a
%major difficulty in many situations, it is impossible for general 
%stochastic systems, especially for queueing
%systems, communication networks, manufacturing systems where systems
%evolve over time in a very complicated manner [e.g. Ho and Cao (1991)]. 
%Recently, there have
%developed two methods, the {\em perturbation analysis} [see, e.g. Ho 
%and Cao (1991)] and the {\em likelihood ratio} method [see, e.g. Reiman 
%and Weiss (1986)], for derivative estimation of the systems that have just
%been mentioned. Although these two methods can yield unbiased estimates 
%of the derivative for a large class of systems,
%they are not general methods. Furthermore, the variance of the
%likelihood ratio method may grow linearly in computation effort. 
On the other hand, although its convergence is slower, the KW algorithm 
requires no detailed information on the function $J(\theta)$. It 
is simple to use and applicable to a wide range of problems. Kesten (1958) 
suggests that the stepsize $a_{n}$ be chosen according to the 
fluctuation in the signs of $g_{n}$ and $h_{n}$. A few of other 
techniques for the acceleration of stochastic approximation algorithms 
can be found in Wasan (1986). None of these 
accelerating techniques can improve the rate of convergence of the
algorithms under study. 

In this paper, we are interested in the acceleration of the KW algorithm 
and the MD algorithm through controlling the estimation of the 
derivative using finite differences.
%Particularly, we examine the convergence of the MD algorithm under the 
%symmetric finite difference (\ref{kw_h}). 
Furthermore, we
consider the employment of the scheme of common random 
numbers (CRN) for improving the convergence of the algorithms --- that is, the random factors 
$\xi_{1,n}$ and $\xi_{2,n}$ are chosen in such a manner that $\xi_{1,n}=
\xi_{2,n}=\xi_{n}$. Implementation of CRN in Monte 
Carlo optimization is rather straightforward. The term ``Monte Carlo optimization'' is used 
here to refer to the procedure of finding the optimal solutions 
through computer simulation where the random factors, represented 
through a sequence of psuedo-random numbers, can be controlled 
[see Bratley et al. (1983)]. Computer simulation 
is often necessary when the form of $L(X(\theta,\xi))$ is too complicated.
This is the case when $L(X(\theta,\xi))$ represents a performance
measure of a stochastic system such as queueing systems, manufacturing
systems, transportation systems, and communications networks [see, e.g.
Ho and Cao (1991);
Bratley et al. (1983); Law and Kelton (1982)]. In Section 6,
we will give an example where the scheme of common random numbers is
feasible.

We use the term CRN in a much narrow sense.
The term CRN has more general and sometimes ill-posed meaning than 
we intend in this paper [see Glasserman and Yao (1992)]. In this paper,
CRN simply refers to that simulation experiments be performed with 
the same stream of random numbers. As far as the KW algorithm (\ref{kw}) 
or the MD algorithm (\ref{md}) is 
concerned, CRN requires that estimates of $J(\theta+\delta)$ 
and $J(\theta-\delta)$ be obtained from simulation experiments 
using the same stream of random numbers $\{\xi_{n}\}$. Let 
$F(\theta,x)$ denote the distribution function of $X(\theta,\xi)$.
Then any experiments with $h_{n}$ constructed in the 
following form conform the CRN requirement:
\begin{equation}
    \frac{L(Y_{1}(\theta_{n},\delta_{n},\xi_{n}))
            -L(Y_{2}(\theta_{n},\delta_{n},\xi_{n}))}{2\delta_{n}} 
 \label{fd_g}
\end{equation}
where the marginal distributions of $Y_{1}(\theta_{n},\delta_{n},\xi_{n})$ 
and $Y_{2}(\theta_{n},\delta_{n},\xi_{n})$ are $F(\theta_{n}+\delta_{n},x)$ 
and $F(\theta_{n}-\delta_{n},x)$, respectively. Note that the joint 
distribution of $Y_{1}(\theta_{n},\delta_{n},\xi_{n})$ and 
$Y_{2}(\theta_{n},\delta_{n},\xi_{n})$ is left open, which may be 
used to improve
the estimation variance. For a distribution function $F(\theta,x)$, its 
inverse function is defined as $F^{-1}(\theta,x) \stackrel{\rm def}{=} 
\inf\{u\;|\;F(\theta,u)> x\}$. Cambanis and Simons (1976) and Whitt (1976)
proved that the variance of (\ref{fd_g}) is minimized when 
$Y_{1}(\theta_{n},\delta_{n},\xi_{n})=F^{-1}(\theta_{n}+
\delta_{n},\xi_{n})$ and $Y_{2}(\theta_{n},\delta_{n},\xi_{n})=
F^{-1}(\theta_{n}-\delta_{n},\xi_{n})$. In this paper we assume 
that the form of $L(X(\theta,\xi))$ is given and fixed. The term 
CRN merely refers to the special choice of $\xi_{1,n}=\xi_{2,n}$. 
We will show that the use of CRN can significantly increase the rate of 
convergence for the KW algorithm (\ref{kw}) or the MD algorithm (\ref{md}) from $n^{-1/3}$ to at least 
$n^{-2/5}$. For a large class of functions, the rate can be increased to $n^{-1/2}$, 
the best possible rate for stochastic approximation algorithms. 
CRN increases the rate of convergence by reducing the variance of 
$h_{n}$. Let $Var[X]$ denote the mathematical variance of a random 
variable $X$. Assume that $Var[L(X(\theta,\xi))]$ is continuous in 
$\theta$, is bounded from below by a positive constant and from above
by a constant. Then if $\xi_{1,n}$ and $\xi_{2,n}$ are 
independent, the variance of $h_{n}$ is
\[   (Var[L(X(\theta_{n}+\delta_{n},\xi_{1,n}))]+
      Var[L(X(\theta_{n}-\delta_{n},\xi_{2,n}))])/(2\delta_{n})^{2}  
     =O(1/\delta_{n}^{2})    \]
which grows quadratically as $\delta_{n}$ goes to zero. We say a variable 
$f(s)=O(s)$ if $|f(s)/s|\leq C, C>0$ is a constant independent of $s$ 
($f(s)=o(s)$ if $\lim|f(s)/s| = 0$ when $s$ goes to zero or infinity 
depending on the context). It is such a large variance of $h_{n}$
that slows down the convergence rate since, when $\delta_{n}$ is 
suitably chosen, the rate would be $n^{-1/2}$ if the variance of 
$h_{n}$ is bounded. As we will show later, the convergence 
rate for (\ref{kw}) depends on how fast the variance of $h_{n}$
goes to infinity. The slower the variance goes to infinity, 
the faster the convergence rate for (\ref{kw}) is. 
CRN 
has been observed effective for variance reduction for decades. It is 
perhaps the most popular method for variance reduction [Bratley et al.
 (1983); Conway (1963); Fishman (1974); Hammersley and 
Handscomb (1964); Heikes et al. (1976); Kleijnen (1974); 
Law and Kelton (1982)]. 
%Heidelberger and Iglehart (1979) first attempted 
%to study theoretically the use of CRN for general stochastic processes. 
%Other related work can be found in Ho and Cao (1991),  Glasserman and Yao (1992); Glynn (1985); 
%Mitchell (1973); Rubinstein and Samorodnitsky (1985);
%Wright and Ramsay (1979). As for the KW algorithm (2), we will show that
%the use of CRN can reduce the variance of $h_{n}$ to $O(1/\delta_{n})$ or
%even $O(1)$. As a result of such a variance reduction, the convergence 
%rate of the KW algorithm is substantially increased.

The rest of the paper is arranged as follows: In Section 2 we examine 
the rates of convergence for the KW algorithm 
under a very general setting that covers many interesting situations. 
the analysis is extended to the MD algorithms in Section 4.
In Section 4 we show that the use of CRN can reduce the variance of
$h_{n}$ by orders of magnitude, which in turn
accelerates the convergence of the KW algorithm. In Section 5, we examine 
the rate of convergence for the MD algorithm under CRN.
In Section 6, we extend the results to multivariates. A practical
example is given to illustrate the feasibility of applying CRN in practice.
Finally, a summary is provided in Section 7.

\vspace{5mm}

{\bf 2. Rates of convergence for the KW algorithm.}\hspace{3mm}
In this section, we examine the rates of convergence for the KW algorithm 
(\ref{kw}) under general assumptions on $h_{n}$. We do not assume that $h_{n}$ 
is of the form (\ref{kw_h}). We will see later that such a treatment covers several 
important cases.

Assume that $\delta_{n}>0$ goes to zero as $n\rightarrow\infty$ and, for
$n\geq n_{0}>1$, $h_{n}$ satisfies the following assumptions:
\begin{equation}
   E[h_{n}|\theta_{n}]=J'(\theta_{n})+\Delta_{n}, \;\;\;
           |\Delta_{n}|\leq b\delta_{n}^{\beta},
\label{h1}
\end{equation}
and
\begin{equation}
     Var[h_{n}|\theta_{n}]\leq c\delta_{n}^{\gamma},
\label{h2}
\end{equation}
where $b, c, \beta$ are real nonnegative numbers, $\gamma\in R$. 
The form of (\ref{h1}) assures that $h_{n}$ is an asymptotically 
unbiased estimate of $J'(\theta)$ when $\beta>0$. When $\gamma>0$, 
the variance of the estimate 
goes to zero as $n\rightarrow\infty$. This is generally impossible 
in practice. When $\gamma =0$ such as in the RM 
algorithm, the variance is bounded. In the case that $\gamma<0$, e.g. 
$\gamma=-2$ if $h_{n}$ is defined by (\ref{kw_h}) and if $\xi_{1,n}$ and 
$\xi_{2,n}$ are independent, the variance of $h_{n}$
goes to infinity. Next we examine the convergence and the rate 
of convergence for the KW algorithm (\ref{kw}). 
The commonly used criterion for measuring the convergence of a stochastic 
sequence $\{\theta_{n}\}$ is the {\em root mean square error} (RMSE) 
defined as
\[    RMSE_{\theta_{n}} = (E[(\theta_{n}-\theta^{*})^{2}])^{1/2}.  \]
If $RMSE_{\theta_{n}} = O(n^{-s}), s>0$, we say that $\{\theta_{n}\}$ 
converges at the rate of $n^{-s}$ or the convergence rate for 
$\{\theta_{n}\}$ is $n^{-s}$.
\vspace{3 mm}

We need the next lemma that was due to Chung (1954) and was formulated
in the present form by Fabian (1971).
\vspace{3 mm}

{\sc Lemma 1}. {\em Let $s, t, B, A_{n}, b_{n}$ be real numbers, $0<s\leq 1$, 
$t\geq 0$, $B>0$. Define $b_{+}=0$ if $s<1$ and $b_{+}=t$ if $s=1$ and 
assume that $c=\lim_{n\rightarrow\infty} A_{n}-b_{+}$ exists and is finite. 
If for $n\geq n_{0}$,
\[ b_{n+1}\leq b_{n}(1-\frac{A_{n}}{n^{s}})
           +\frac{B}{n^{s+t}}   \]
and if $c>0$, then
\[    \lim_{n\rightarrow \infty}\sup n^{t}b_{n} \leq B/c.     \]
The statement remains valid if all the inequalities are reversed and 
$\lim\sup$ is replaced by $\lim\inf$}.
\vspace{3 mm}

The following Theorems 1 and 2 give the convergence rate for the KW 
algorithm (\ref{kw}) with $h_{n}$ satisfying (\ref{h1})-(\ref{h2}):
\vspace{3 mm}

{\sc Theorem 1}. {\em Assume that $\{\theta_{n}\}$ is determined by (\ref{kw}) and
\begin{itemize}
\item[(A1).] $a_{n}=a n^{-\alpha}, \delta_{n}=d n^{-\eta}$, $0<\alpha
            \leq 1$, $\eta>0$, $a, d > 0$;
\item[(A2).] $J(\theta)$ is increasing for $\theta< 
            \theta^{*}$ and decreasing for $\theta>\theta^{*}$, and there 
            exist two constants $K_{1}, K_{2}$, $0<K_{1}\leq K_{2}
            <\infty$, such that for all $\theta\in \Theta$,
            \[ K_{1}|\theta-\theta^{*}|\leq |J'(\theta)|\leq
               K_{2}|\theta-\theta^{*}|;   \]
\item[(A3).] conditioned on $\theta_{n}$, $h_{n}$ at the $n$th iteration is
             independent of those at the other iterations.
\end{itemize}
Then, if $\sigma=(1/2)\min\{\alpha+\gamma\eta,2\beta\eta\}$ and $0< 
\sigma < aK_{1}$, we have
\begin{equation}
    \lim_{n\rightarrow\infty}\sup 
             n^{2\sigma}E[(\theta_{n}-\theta^{*})^{2}]\leq C 
\label{kw_rate1}
\end{equation}
where $C > 0$ is a constant. The convergence rate for $RMSE_{\theta_{n}}$
is at least $n^{-\sigma}$}.

{\sc Proof}. Without loss of generality, we assume that $\theta^{*}=0$.
Then 
\begin{eqnarray*}
 E[\theta_{n+1}^{2}] & = &E[\theta_{n}^{2}]-2a_{n}E[\theta_{n}h_{n}]
            +a_{n}^{2}E[h_{n}^{2}]          \nonumber \\
 & = & E[\theta_{n}^{2}]-2a_{n}E[\theta_{n}(J'(\theta_{n})
            +\Delta_{n})]+a_{n}^{2}((E[h_{n}])^{2}+Var[h_{n}]).
\end{eqnarray*}
According to (\ref{h1})-(\ref{h2}), we have
\begin{equation}
  E[\theta_{n+1}^{2}]  \leq 
     E[\theta_{n}^{2}]-2a_{n}E[\theta_{n}J'(\theta_{n})]
      +2b a_{n}\delta_{n}^{\beta}E[|\theta_{n}|] 
  +2a_{n}^{2}(E[J'(\theta_{n})]^{2}         
     +b^{2}\delta_{n}^{2\beta})+c a_{n}^{2}\delta_{n}^{\gamma}. 
\label{bound1}
\end{equation}
By Assumption (A2), $\theta_{n}J'(\theta_{n}) \geq 0$ and
\begin{equation}
           \theta_{n}J'(\theta_{n})\geq K_{1}\theta_{n}^{2}, \;\; 
             (J'(\theta_{n}))^{2}\leq K_{2}^{2}\theta_{n}^{2}.
\label{bound_J}
\end{equation}
Furthermore, for any $\epsilon_{n} > 0$, $|\theta_{n}|\leq \epsilon_{n}+
\theta_{n}^{2}/\epsilon_{n}$ and consequently
\[  E[|\theta_{n}|] \leq 
      \epsilon_{n}+\frac{1}{\epsilon_{n}}E[\theta_{n}^{2}].  \]
By setting $0< \epsilon< 1$ and
\[   \epsilon_{n} =  \frac{2b\delta_{n}^{\beta}}{K_{1}\epsilon}, \]
we have
\begin{equation}  E[|\theta_{n}|] \leq 
       \frac{2b\delta_{n}^{\beta}}{K_{1}\epsilon}
          +\frac{K_{1}\epsilon}{2b\delta_{n}^{\beta}}E[\theta_{n}^{2}].
\label{bound_th}
\end{equation}
Substituting (\ref{bound_J}) and (\ref{bound_th}) into (\ref{bound1}), we obtain 
\begin{equation}  
  E[\theta_{n+1}^{2}]\leq E[\theta_{n}^{2}][1-(2-\epsilon)K_{1}a_{n}
         +2K_{2}^{2}a_{n}^{2}]+2b^{2}a_{n}^{2}\delta_{n}^{2\beta}
      +c a_{n}^{2}\delta_{n}^{\gamma}
    +\frac{4b^{2}}{K_{1}\epsilon}a_{n}\delta_{n}^{2\beta}.
\label{bound_th2}
\end{equation}
According to (\ref{bound_th2}), also noting Assumption (A1), we can choose an 
$n_{1}\geq n_{0}>1$ such that for all $n\geq n_{1}$
\[   E[\theta_{n+1}^{2}]\leq E[\theta_{n}^{2}](1-\frac{A_{n}}{n^{\alpha}})
         +\frac{B}{n^{\alpha+2\sigma}}      \]
where
\[   A_{n}= (2-\epsilon)a K_{1} 
             -\frac{2K_{2}^{2}a^{2}}{n^{\alpha}},\;\;\;
     B=2 a^{2}b^{2}d^{2\beta}+c a^{2}d^{\gamma}
           +\frac{4a b^{2}d^{2\beta}}{K_{1}\epsilon}. \]
If $aK_{1}>\sigma$, we can always choose $\epsilon>0$ so small that
$(2-\epsilon)aK_{1}>2\sigma$. Applying Lemma 1, we obtain (\ref{kw_rate1}) with $C=
B/((2-\epsilon)aK_{1})$ if $\alpha<1$ and $C=B/((2-\epsilon)aK_{1}-2\sigma)$
if $\alpha=1$.
\hspace{0.1in}\rule{1.5mm}{3.0mm}
\vspace{3 mm}

It follows directly from Theorem 1 that $\{\theta_{n}\}$ converges to 
$\theta^{*}$ as long as $\sigma > 0$, or equivalently, as long as $\alpha
+\gamma\eta > 0$. When $\alpha+\gamma\eta \leq 0$ which is possible only 
when $\gamma < 0$, the variance of $h_{n}$ grows 
to infinity at the rate of $n^{t}$ with $t=-\gamma\eta \geq \alpha$. It is 
obvious from (\ref{kw}) that $\{\theta_{n}\}$ does not converge.
Another extreme case is that $\gamma > 0$. In this case, $\sigma$ can be
made arbitrarily large by choosing appropriate $\eta$.
The convergence rate for $\{\theta_{n}\}$ can be made arbitrarily large
if $\eta$ can take any value. In fact, by setting $\eta\rightarrow \infty$ 
in (\ref{bound_th2}) (or equivalently, $\delta_n\rightarrow 0$) and
 $a_{n}=a$ such that $ 0 < q=1-(2-\epsilon)K_{2}a
+2K_{2}^{2}a^{2} < 1$ for sufficiently large $n$, we have
\[  E[\theta_{n+1}^{2}] \leq q E[\theta_{n}^{2}].  \]
The convergence rate for the sequence $\{\theta_{n}\}$ is that of a 
geometric progression. Unfortunately, this is
a very special case. One should not expect $\gamma > 0$ in practice.
Both of the situations $\gamma > 0$ and $\alpha+\gamma\eta\leq 0$ 
are too special to deserve further study. The most interesting case 
is when $\gamma$ satisfies $-\alpha/\eta < \gamma \leq 0$.
\vspace{3 mm}

Theorem 1 shows that, when $h_{n}$ satisfies (\ref{h1})-(\ref{h2}), $\{\theta_{n}\}$ 
converges with probability one to the optimal parameter $\theta^{*}$ 
at a rate of at least $n^{-\sigma}$. We can further prove that 
$\{\theta_{n}\}$ converges exactly at this rate as interpreted in the 
following Theorem 2.
\vspace{3 mm}

{\sc Theorem 2}. {\em Assume that Assumptions (A1)-(A3) are satisfied.
\begin{enumerate}
\item If $\alpha+\gamma\eta < 2\beta\eta$, $aK_{2}>\sigma$, 
      $E[h_{n}|\theta_{n}]=J'(\theta_{n})+\Delta_{n}$, $|\Delta_{n}|\leq
      b\delta_{n}^{\beta}$, and $Var[h_{n}|\theta_{n}] \geq 
      c\delta_{n}^{\gamma}$, we have
      \[    \lim_{n\rightarrow\infty}\inf 
             n^{2\sigma}E[(\theta_{n}-\theta^{*})^{2}]\geq C_{1}   \]
       where $C_{1}>0$ is a constant.
\item If $\alpha+\gamma\eta \geq 2\beta\eta$, $E[h_{n}|\theta_{n}]=
      J'(\theta_{n})+b\delta_{n}^{\beta}(1+\varepsilon_{n})$,
      and $J'(\theta_{n})=(\theta_{n}-\theta^{*})(K_{3}+\tau_{n})$, 
      $\varepsilon_{n}=o(1)$ and $\tau_{n}=o(1)$ uniformly 
      as $n\rightarrow\infty$, $K_{3}=J''(\theta^{*})>0$, 
      $\sigma < aK_{3}$, then
      \[   \lim_{n\rightarrow\infty}\sup 
            n^{\sigma}E[\theta_{n}-\theta^{*}]\leq - C_{2}      \]
      where $C_{2}>0$ is a constant.
\end{enumerate}}

{\sc Proof}. Let consider the first statement. For simplicity and without loss
of generality, we assume
$\theta^{*}=0$. Parallel to the derivation of (\ref{bound1}) we have
\begin{eqnarray*}
 E[\theta_{n+1}^{2}] 
  =  E[\theta_{n}^{2}]-2a_{n}E[\theta_{n}(J'(\theta_{n})
          +\Delta_{n})] +a_{n}^{2}\{(E[h_{n}])^{2}
                 +Var[h_{n}]\} 
\end{eqnarray*}
which implies
%\begin{equation}
\[ E[\theta_{n+1}^{2}] 
  \geq E[\theta_{n}^{2}]-2a_{n}E[\theta_{n}J'(\theta_{n})]
           -2b a_{n}\delta_{n}^{\beta}E[|\theta_{n}|]
          +c a_{n}^{2}\delta_{n}^{\gamma}. 
\]%\label{bound_th2}
%\end{equation}
Assumption (A2) implies that $0\leq\theta_{n}J'(\theta_{n}) \leq 
K_{2}\theta_{n}^{2}$ which, together with (\ref{bound_th}) where $K_{1}$ is replaced 
by $K_{2}$, shows that 
\begin{equation}
   E[\theta_{n+1}^{2}]\geq E[\theta_{n}^{2}](1-(2-\epsilon)K_{2}a_{n})
        +c a_{n}^{2}\delta_{n}^{\gamma}
        -\frac{4b^{2}}{K_{2}\epsilon}a_{n}\delta_{n}^{2\beta}.
\label{bound_th3}
\end{equation}

If $\alpha+\gamma\eta < 2\beta\eta$, there
exists an $n_{0}>1$ such that when $n\geq n_{0}$
\[    c a_{n}^{2}\delta_{n}^{\gamma}-\frac{4b^{2}}{K_{2}\epsilon}
      a_{n}\delta_{n}^{2\beta}
    \geq \frac{1}{2} c a_{n}^{2}\delta_{n}^{\gamma}.  \]
Therefore, we know from (\ref{bound_th3}) that when $n\geq n_{0}$
\[   E[\theta_{n+1}^{2}]\geq E[\theta_{n}^{2}](1-\frac{A_{n}}{n^{\alpha}})
         +\frac{ca^{2}d^{\gamma}}{2n^{\alpha+2\sigma}}.      \]
Since $aK_{2}>\sigma$, we can always choose $\epsilon>0$ so small that
$(2-\epsilon)aK_{2}>2\sigma$. The first statement of the theorem follows
from applying Lemma 1 with $C_{1}= ca^{2}d^{\gamma}/(2(2-\epsilon)aK_{2})$ 
if $\alpha<1$ and $C_{1}= ca^{2}d^{\gamma}/(2(2-\epsilon)aK_{2}-4\sigma)$ if 
$\alpha=1$.

If $\alpha+\gamma\eta \geq 2\beta\eta$, we know that $\sigma=\beta\eta>0$
and
\begin{eqnarray}
     E[\theta_{n+1}] &= & E[\theta_{n}] - a_{n}E[J'(\theta_{n})]-
               b a_{n}\delta_{n}^{\beta}(1+\varepsilon_{n}) \nonumber \\
   &=& E[\theta_{n}](1-K_{3}a_{n})+a_{n}E[\tau_{n}\theta_{n}]-
               b a_{n}\delta_{n}^{\beta}(1+E[\varepsilon_{n}])  
\label{eth1}
\end{eqnarray}
Define $z_{n}=n^{\sigma}E[\theta_{n}]$. Then (\ref{eth1}) shows that
\begin{eqnarray*}
   z_{n+1} & = & n^{\sigma}(1+\frac{1}{n})^{n^{\sigma}}
        [E[\theta_{n}](1-K_{3}a_{n})+a_{n}E[\tau_{n}\theta_{n}]
        -b a_{n}\delta_{n}^{\beta}(1+\varepsilon_{n})] \nonumber \\
  & = & z_{n}[1+\frac{\sigma}{n}-K_{3}a_{n}
        -\frac{\sigma}{n}K_{3}a_{n}+O(\frac{1}{n^{2}})] 
   +(1+\frac{1}{n})^{\sigma} n^{\sigma}(a_{n}E[\tau_{n}\theta_{n}]
       -ba_{n}\delta_{n}^{\beta}(1+E[\varepsilon_{n}])).
\end{eqnarray*}
Denote 
\[    A_{n} = 1+\frac{\sigma}{n}-K_{3}a_{n}
         -\frac{\sigma}{n}K_{3}a_{n}+O(\frac{1}{n^{2}}),   \]
\[    B_{n} = (1+\frac{1}{n})^{\sigma} n^{\sigma}(ba_{n}
         \delta_{n}^{\beta}(1+E[\varepsilon_{n}])-a_{n}E[\tau_{n}\theta_{n}]).\]
Then
\begin{equation}
       z_{n+1}= A_{n}z_{n}-B_{n}.
\label{z1}
\end{equation}
Note that $a_{n}=an^{-\alpha}, 0<\alpha\leq 1, \delta_{n}=dn^{-\eta}, 
aK_{3}>\sigma$. We may choose $\tilde{A}_{1}, \tilde{A}_{2}>0, n_{1}>1$ 
such that, for all $n \geq n_{1}$,
\begin{equation}
   0\leq 1-\frac{\tilde{A}_{1}}{n^{\alpha}} 
       \leq A_{n}=1+\frac{\sigma}{n}-\frac{aK_{3}}{n^{\alpha}}
     -\frac{a\sigma K_{3}}{n^{1+\alpha}}+O(\frac{1}{n^{2}})
     \leq 1-\frac{\tilde{A}_{2}}{n^{\alpha}}.
\label{A1}
\end{equation}
Since Assumptions (A1)-(A3) in Theorem 1 are satisfied, 
$\lim_{n\rightarrow\infty}\sup n^{2\sigma}E[\theta_{n}^{2}]\leq C$ which
implies that 
\[  \lim_{n\rightarrow\infty}\sup n^{\sigma}|E[\theta_{n}]|
   \leq \lim_{n\rightarrow\infty}\sup (n^{2\sigma}E[\theta_{n}^{2}])^{1/2}
   \leq \sqrt{C}.  \]
According to the assumptions that $\epsilon_{n}=o(1)$, $\tau_{n}=o(1)$ 
uniformly as $n\rightarrow\infty$, and $\delta_{n}^{\beta}=
d^{\beta}n^{-\sigma}$. There exists an $n_{2}> 1$ 
such that, when $n\geq n_{2}$, 
\begin{eqnarray}
   B_{n}& = &(1+\frac{1}{n})^{\sigma} n^{\sigma}(ba_{n}\delta_{n}^{\beta}(1
      +E[\varepsilon_{n}])-a_{n}E[\tau_{n}\theta_{n}])     \nonumber \\
   & \geq &(1+\frac{1}{n})^{\sigma} n^{\sigma}(ba_{n}\delta_{n}^{\beta}(1
      +E[\varepsilon_{n}])-a_{n}E[|\tau_{n}|]E[|\theta_{n}|])     \nonumber \\
   & \geq & (1+\frac{1}{n})^{\sigma} n^{\sigma}
             \frac{1}{2}ba_{n}\delta_{n}^{\beta}
         \geq \frac{1}{2}n^{\sigma}ba_{n}\delta_{n}^{\beta}.    
\label{Bn1}
\end{eqnarray}
Let $n_{0}=\max\{n_{1},n_{2}\}$. Then, from (\ref{z1}) we know that for all 
$n\geq n_{0}$
\begin{equation}
    z_{n} = z_{n_{0}}\prod^{n}_{i=n_{0}}A_{i}
            -\sum_{i=n_{0}}^{n-1}B_{i}\prod_{j=i+1}^{n}A_{j}-B_{n}.
\label{z2}
\end{equation}
Since $0<\alpha\leq 1$, (\ref{A1}) shows that
\begin{equation}  0\leq \lim_{n\rightarrow\infty} \prod^{n}_{i=n_{0}}A_{i}
     \leq  \lim_{n\rightarrow\infty}
     \prod^{n}_{i=n_{0}}(1-\frac{\tilde{A}_{2}}{i^{\alpha}}) = 0.   
\label{limA}
\end{equation}
Furthermore, $\lim_{n\rightarrow\infty}B_{n} = 0$, and (\ref{A1}) and (\ref{Bn1}) imply
that
\begin{eqnarray}
     \sum_{i=n_{0}}^{n-1}B_{i}\prod_{j=i+1}^{n}A_{j}
      \geq \sum_{i=n_{0}}^{n-1}\frac{abd^{\beta}}{2i^{\alpha}}
              \prod_{j=i+1}^{n}A_{j}
      \geq \frac{abd^{\beta}}{2n^{\alpha}}\sum_{i=n_{0}}^{n-1}A_{n}^{n-i}
      \geq \frac{abd^{\beta}}{2n^{\alpha}}\sum_{i=n_{0}}^{n-1}
        (1-\frac{\tilde{A}_{1}}{n^{\alpha}})^{n-i}.
\label{BA}
\end{eqnarray}
On the other hand,
\[   \lim_{n\rightarrow\infty}\frac{1}{n^{\alpha}}\sum_{i=n_{0}}^{n-1}
        (1-\frac{\tilde{A}_{1}}{n^{\alpha}})^{n-i}
      =\left\{ \begin{array}{ll}
           1, & \mbox{if }0<\alpha<1 \\ 
           1-e^{-\tilde{A}_{1}}, & \mbox{if }\alpha=1.
       \end{array} \right.  \]
Substituting the preceding inequality, (\ref{limA}) and (\ref{BA}) into (\ref{z2}), we see
that
\[     \lim_{n\rightarrow\infty}\sup z_{n} \leq - C_{2}   \]
with $C_{2}=(1/2)abd^{\beta}(1-e^{-\tilde{A}_{1}})>0$. This is exactly what we want to 
prove.
\hspace{0.1in}\rule{1.5mm}{3.0mm}
\vspace{3 mm}

Theorems 1 and 2 show that the convergence rate for $\{\theta_{n}\}$ is
generally $n^{-\sigma}$. If we are free to choose the positive numbers 
$\alpha, \eta$, it follows directly from Theorem 1 that
\vspace{3 mm}

{\sc Corollary 1}. {\em Assume that $h_{n}$ satisfies (\ref{h1})-(\ref{h2}) and $\gamma
\leq 0$. Under Assumptions (A2)-(A3) in Theorem 1, the best possible
convergence rate for the KW algorithm (\ref{kw}) is $n^{-\beta/(2\beta-\gamma)}$ 
which is achieved by setting $\alpha = 1$, $\eta= 1/(2\beta-\gamma)$, 
and by choosing appropriate $a, d>0$.}
\vspace{3 mm}

For the KW algorithm (\ref{kw}) with $h_{n}$ defined by (\ref{kw_h}), assume that 
$J(\theta)$ is continuously differentiable of order up to three and the 
third order derivative $J'''(\theta)$ is uniformly bounded on $\Theta$, we have
\[   E[h_{n}|\theta_{n}] = \frac{J(\theta_{n}+\delta_{n}) -J(\theta_{n}-
      \delta_{n})}{2\delta_{n}} = J'(\theta_{n})+
     \frac{1}{6}J'''(\tilde{\theta}_{n})\delta_{n}^{2} 
             = J'(\theta_{n})+O(\delta_{n}^{2})   \]
where $\tilde{\theta}_{n}\in[\theta_{n}-\delta_{n},\theta_{n}+\delta_{n}]$. 
In this case, $\beta=2$. If the assumptions (\ref{h2}) and (A1)-(A3) are 
satisfied and if the positive number $a$ is chosen sufficiently large, 
we know from Theorems 1 and 2 that
\begin{equation}  
   \sigma = \frac{1}{2}\min\{\alpha+\gamma\eta,4\eta\}.
\label{sigma1}
\end{equation}
%\vspace{3 mm}

If we use the one-sided finite-difference approximation in (\ref{kw}):
\begin{equation}
    h_{n}=\frac{L(X(\theta_{n}+\delta_{n},\xi_{1,n}))
                   -L(X(\theta_{n},\xi_{2,n}))}{\delta_{2,n}}
\label{kw_1sideh}
\end{equation}
and if $J(\theta)$ is twice continuously differentiable and the second order
derivative $J''(\theta)$ is bounded on $\Theta$, then for any $\theta_{n},
\delta_{n}$ there exists a $\hat{\theta}_{n}\in[\theta_{n},\theta_{n}
+\delta_{n}]$ such that
\[   E[h_{n}|\theta]=\frac{J(\theta_{n}+\delta_{n})-J(\theta_{n})}{\delta_{n}} 
       = J'(\theta_{n})+\frac{1}{2}J''(\hat{\theta}_{n})\delta_{n} 
            = J'(\theta_{n})+O(\delta_{n}). \]
Therefore, $\beta=1$. Under the same conditions as those in the previous 
case we know that
\begin{equation}
    \sigma = \frac{1}{2}\min\{\alpha+\gamma\eta,2\eta\}.
\label{sigma2}
\end{equation}
It is clear from (\ref{sigma1}) and (\ref{sigma2}) that, under the same condition for the
variance $Var[h_{n}|\theta_{n}]$, the convergence rate of the KW algorithm
is faster when symmetric differences are used than that when 
one-sided differences are used. Corollary 1 shows that the best possible 
convergence rate depends on two factors---how fast the bias decreases 
to zero and how slow the variance
increases to infinity. Using symmetric finite difference
(\ref{kw_h}) instead of the one-sided finite difference (\ref{kw_1sideh}) can 
 reduce the bias of $h_{n}$. 
%In Section 4 we will demonstrate 
%how to use CRN to reduce the variance of $h_{n}$, which can drastically 
%increase the convergence rate for the KW algorithm. 
%Combining 
%the preceding remarks with Corollary 1, 
To summarize, we have the following 
conclusion which will be used later.
\vspace{3mm}

{\sc Corollary 2}. {\em Suppose that (A1)-(A3) are satisfied. If
\begin{itemize}
\item[(A4).] $J(\theta)$ is continuously differentiable of order up to 
             three and the third order derivative $J'''(\theta)$ is 
             bounded on $\Theta$,
\end{itemize}
then the best possible convergence rate for the KW algorithm (\ref{kw}) with 
$h_{n}$ defined in (\ref{kw_h}) is $n^{-2/(4-\gamma)}$. If
\begin{itemize}
\item[(A5).] $J(\theta)$ is twice continuously differentiable and the 
             second order derivative $J''(\theta)$ is bounded on $\Theta$,
\end{itemize}
then the best possible convergence rate for the KW algorithm (\ref{kw}) with 
$h_{n}$ defined in (\ref{kw_1sideh}) is $n^{-1/(2-\gamma)}$.}

\vspace{5mm}

{\bf 3. Rates of convergence for the MD algorithm.}\hspace{3mm}
The rate of convergence of the MD algorithms was established 
by Nemivoski et al. (2009) when the $h_n$ in (\ref{md}) is an unbiased
estimate of the derivative, and by Duchi et al. (2012, 2013) when the $h_n$ is 
approximated by the one-sided finite difference (\ref{kw_1sideh}). 
In this section, we examine the rate of convergence of the MD algorithm 
for general $h_n$. Again, we only assume that $h_n$ satisfies (\ref{h1})-(\ref{h2}). 
For notational consistence, 
the norm $||.||$ in (\ref{md}) is taken as the $l_2$ norm. Its dual norm 
$||x||_{*}:=\sup_{||y||\leq 1}y^Tx$ is also the $l_2$ norm.
Define
\[ \hat{\theta}_n = \frac{1}{n}\sum_{i=1}^n \theta_i. \]
We next examine the convergence of $J( \hat{\theta}_n)$.
\vspace{3 mm}

{\sc Theorem 3}. {\em Assume that $\{\theta_{n}\}$ is determined by (\ref{md}), and
\begin{itemize}
\item[(B1).] $\psi(\theta)$ is strongly convex, $\Theta$ is compact and convex,
and there exists $r>0$ such that $D(\theta^{*},\theta)\leq (1/2)r^2, r > 0$ for all $\theta\in\Theta$;
\item[(B2).] $L(X)$ is closed convex, and there 
            exist two constants $K_{1}, K_{2}, 0<K_{1}\leq K_{2}
            <\infty$, such that for all $\theta\in \Theta$,
            \[ K_{1}|\theta-\theta^{*}|\leq |J'(\theta)|\leq
               K_{2}|\theta-\theta^{*}|;   \]
\item[(B3).] conditioned on $\theta_{n}$, $h_{n}$ at the $n$th iteration is
             independent of those at the other iterations.
\end{itemize}
Then
\begin{equation}
    E[J(\hat{\theta}_{n})-J(\theta^{*})] \leq \frac{C_1}{n a_n}+\frac{C_2}{n}\sum_{i=1}^n a_i \delta_i^\gamma
+\frac{C_3}{n}\sum_{i=1}^n a_i 
+\frac{C_4}{n}\sum_{i=1}^n a_i \delta_i^{2\beta}
+ \frac{C_5}{n}\sum_{i=1}^n \delta_i^\beta,
\label{eqn_master}
\end{equation}
where
\[ C_1 =  \frac{r^2}{2},\;\; C_2 = \frac{c}{2\kappa}, \;\; C_3 =\frac{K_2^2r^2}{2\kappa^2},\;\; C_4=\frac{b^2}{\kappa},\;\; C_5= \frac{br}{\sqrt{2\kappa}}.\]
}

{\sc Proof}. Under Assumptions (B1)-(B3), we know from Duchi et al. (2013), eqn. (13), that
\begin{equation}
 J(\hat{\theta}_n) - J(\theta^{*}) \leq \frac{r^2}{2n a_n}+\frac{1}{2n\kappa}\sum_{i=1}^n a_i h_i^2 - \frac{1}{n}\sum_{i=1}^n \Delta_i(\theta_i-\theta^{*}).
\label{eqn1}
\end{equation}
Therefore,
\begin{equation}
 E[ J(\hat{\theta}_n) - J(\theta^{*})] \leq \frac{r^2}{2n a_n}+\frac{1}{2n\kappa}\sum_{i=1}^n a_i E[h_i^2] + \frac{1}{n}\sum_{i=1}^n E[|\Delta_i(\theta_i-\theta^{*})|].
\label{eqn2}
\end{equation}
The assumptions (\ref{h1})-(\ref{h2}) give
\begin{equation}
E[h_i^2] = Var[h_i] +(E[J'(\theta_i)+\Delta_i])^2 \leq  c\delta_i^\gamma+2(E[J'(\theta_i)])^2+2b^2\delta_i^{2\beta}
\label{eqn_h1}
\end{equation}
On the other hand, according to Assumptions (B1)-(B2),
\begin{equation}
 (E[J'(\theta_i)])^2 \leq (K_2E[|\theta_i-\theta^{*}|])^2 \leq \frac{K_2^2}{2\kappa}r^2
\label{eqn_h2}
\end{equation}
and
\begin{equation}
 E[|\Delta_i||\theta_i-\theta^{*}|] \leq \frac{ b\delta_i^\beta r}{\sqrt{2\kappa}}.
\label{eqn_h3}
\end{equation}
By combining  (\ref{eqn_h1})-(\ref{eqn_h3}) with (\ref{eqn2}), we obtain
\[
E[ J(\hat{\theta}_n) - J(\theta^{*})]  \leq \frac{r^2}{2n a_n}+\frac{c}{2n\kappa}\sum_{i=1}^n a_i \delta_i^\gamma
+\frac{K_2^2r^2}{2n\kappa^2}\sum_{i=1}^n a_i 
+\frac{b^2}{n\kappa}\sum_{i=1}^n a_i \delta_i^{2\beta}
+ \frac{br}{n\sqrt{2\kappa}}\sum_{i=1}^n \delta_i^\beta,
\]
which is exactly (\ref{eqn_master}). 
\hspace{0.1in}\rule{1.5mm}{3.0mm}
\vspace{3 mm}

A special case of interest is $\gamma=0$, which corresponds to bounded variance of 
derivative estimation. The following  Corollary 3 provides a bound on the convergence of the 
MD algorithm for this case with properly chosen  $\{a_n\}$, $\{\delta_n\}$. 
\vspace{3 mm}

{\sc Corollary 3}. {\em Suppose that (A4), (B1)-(B3) are satisfied. Let $\gamma=0$,
$a_n= a n^{-1/2}$, $\delta_n = d n^{-1}$, $a>0, d>0$.

\noindent (a) For $h_n$ defined in (\ref{kw_1sideh}),
\begin{equation}
 E[J(\hat{\theta}_{n})-J(\theta^{*})] \leq (C_1+2C_2+2C_3)\frac{\max(a,a^{-1})}{\sqrt{n}}+(2.5C_4 d^{2})\frac{a}{n}+(C_5d)\frac{1+\log n}{n}.
\label{md_case0}
\end{equation}

\noindent (b) For $h_n$ defined in (\ref{kw_h}), 
\begin{equation}
 E[J(\hat{\theta}_{n})-J(\theta^{*})] \leq (C_1+2C_2+2C_3)\frac{\max(a,a^{-1})}{\sqrt{n}}+(9C_4 d^{4}/7)\frac{a}{n}+(2C_5 d^{2})\frac{1}{n}.
\label{md_case1}
\end{equation}

}

{\sc Proof}. For  $\gamma=0$, $a_n= a n^{-1/2}$, $\delta_n = d n^{-1}$,
combining the first three terms in (\ref{eqn_master}) gives the first term in (\ref{md_case0}). 
For $h_n$ defined by the one-sided finite difference (\ref{kw_1sideh}), under Assumption (A4), 
we have $\beta=1$. Consequently, the fourth term in  (\ref{eqn_master}) is
\[ \frac{C_4}{n}\sum_{i=1}^n a_i \delta_i^{2} = \frac{C_4 a d^{2}}{n}\sum_{i=1}^n i^{-2.5}
       \leq  \frac{C_4 a d^{2}2.5}{n}.
\]
The last term is
\[  \frac{C_5}{n}\sum_{i=1}^n \delta_i =  \frac{C_5d}{n}\sum_{i=1}^n i^{-1} \leq \frac{C_5d (1+\log n)}{n}.
\]
Combing the previous two inequalities with (\ref{eqn_master}) gives (\ref{md_case0}).

For $h_n$ defined by the symmetric finite difference (\ref{kw_h}), under Assumption (A4), 
we have $\beta=2$. Consequently, the fourth term in  (\ref{eqn_master}) is
\[ \frac{C_4}{n}\sum_{i=1}^n a_i \delta_i^{4} = \frac{C_4 a d^{4}}{n}\sum_{i=1}^n i^{-4.5}
       \leq  \frac{C_4 a d^{4}(9/7)}{n}.
\]
The last term is
\[  \frac{C_5}{n}\sum_{i=1}^n \delta_i^2 =  \frac{C_5d^2}{n}\sum_{i=1}^n i^{-1} \leq \frac{2C_5d^2}{n}.
\]
Combing the previous two inequalities with (\ref{eqn_master}) gives (\ref{md_case1}).
\hspace{0.1in}\rule{1.5mm}{3.0mm}
\vspace{3 mm}

Duchi et al. (2013) investigated the convergence of the MD algorithm using the one-sided finite difference (\ref{kw_1sideh}) 
as an approximation to the derivative. The bound (\ref{md_case0}) is technically the same
as that in Duchi et al. (2013). When the symmetric finite difference (\ref{kw_h}) is used, the $\log n$ factor 
disappears in the last term of (\ref{md_case1}), which indicates that the symmetric 
finite-difference approximation (\ref{kw_h}) leads to a tighter bound
under similar assumptions, which is due to that the symmetric finite difference 
(\ref{kw_h}) typically provides more
accurate estimate of the mean of the derivative than the one-sided ones do. 
Note that Duchi et al. (2012, 2013) implicitly assumes that CRN is used in calculating
the finite difference (\ref{kw_h}) or (\ref{kw_1sideh}) that will be covered in Sections 4 and 5.

It's worth of noting that the rate of convergence for the MD algorithm, as given by 
(\ref{eqn_master}), is $n^{-1/2}$ which is not affected by the choice of finite-difference
approximation, either symmetric or one-sided, to the derivative. This is by design since the MD algorithm 
was originally proposed  for improving the robustness in the choice of stepsizes at the cost of slower convergence.
\vspace{3 mm}

When a finite difference is used to approximate the derivative, it is desirable to have 
$\delta_n\rightarrow 0$ as $n\rightarrow\infty$ to ensure asymptotically 
unbiased estimate of the derivative. In this case, it is possible (and likely in practice!) that the variance of
the estimate goes to infinity. This is a special case of (\ref{eqn_master})
with $\gamma<0$. Therefore, Theorem 3 allows flexibility to cover general cases.

A special situation is when $L(X(\theta_{n}+\delta_{n},\xi_{1,n}))$ and $L(X(\theta_{n}-\delta_{n},\xi_{2,n}))$ 
or  $L(X(\theta_{n},\xi_{2,n}))$ in (\ref{kw_h}) or (\ref{kw_1sideh})
 are sampled independently. In this case, $\gamma=-2$. Assume further that 
$\{a_n\}$ and $\{\delta_n\}$ are specified as in Assumption (A1). Then the right hand side 
of (\ref{eqn_master}) becomes
\[   H(n) := \frac{C_1}{n a_n}+\frac{C_2}{n}\sum_{i=1}^n a_i \delta_i^\gamma
+\frac{C_3}{n}\sum_{i=1}^n a_i 
+\frac{C_4}{n}\sum_{i=1}^n a_i \delta_i^{2\beta}
+ \frac{C_5}{n}\sum_{i=1}^n \delta_i^\beta
 \]
\[ = \frac{C_1}{a n^{1-\alpha}}+\frac{C_2}{n}\sum_{i=1}^n a\delta^{-2} i^{-\alpha+2\eta}
+\frac{C_3}{n}\sum_{i=1}^n a i^{-\alpha}
+\frac{C_4}{n}\sum_{i=1}^n a \delta^{2\beta} i^{-\alpha-\beta\eta}
+ \frac{C_5}{n}\sum_{i=1}^n \delta  i^{-\beta\eta} 
\]
\[ = O(n^{-1+\alpha})+O( n^{-\alpha+2\eta})+O( n^{-\alpha})+O( n^{-\alpha-\beta\eta})+O(n^{-\beta\eta}) 
\]
\[ = O(n^{-\sigma}), \]
where
\[  \sigma  = \min\{ 1-\alpha, \alpha-2\eta, \alpha, \alpha+2\beta \eta, \beta\eta\}
= \min\{ 1-\alpha, \alpha+2\eta, \beta\eta\}. \]         
 For the one-sided finite difference (\ref{kw_1sideh}), $\beta=1$. Then
 \[  \sigma = \min\{ 1-\alpha, \alpha+2\eta, \eta\} \leq 1/4. \]
For the symmetric finite difference (\ref{kw_h}), $\beta=2$. Then
 \[  \sigma = \min\{ 1-\alpha, \alpha+2\eta, 2\eta\} \leq 1/3. \]
 The previous discussion can be summarized in the following Corollary 4.          
\vspace{3 mm}

{\sc Corollary 4}. {\em Assume that Assumptions (A1), (A4), (B1)-(B3) are satisfied, and 
that $L(X(\theta_{n}+\delta_{n},\xi_{1,n}))$ and $L(X(\theta_{n}-\delta_{n},\xi_{2,n}))$ in 
(\ref{kw_h}) (or $L(X(\theta_{n},\xi_{2,n}))$ in (\ref{kw_1sideh}))
are independent. Then
\begin{itemize}
\item[(i)] the best possible rate of convergence for the upper bound $H(n)$ is $n^{-1/4}$
when the one-sided finite difference (\ref{kw_1sideh}) is used,
\item[(ii)] the best possible rate of convergence for the upper bound $H(n)$ is $n^{-1/3}$
when the symmetric finite difference (\ref{kw_h}) is used.
\end{itemize}
}

Note that the rates of convergence are only upper bounds of  $ E[J(\hat{\theta}_{n})]$.
Such rates of convergence are consistent with those for $\{\theta_n\}$.

\vspace{5mm}

{\bf 4. The KW algorithm with CRN.}\hspace{3mm}
In this section, we will show how CRN can
accelerate the convergence of the KW algorithm. For 
 clarity and without getting trapped into unnecessary
tediousness of details, we focus our attention on the case in which
$\xi\in R$ is a real one-dimensional random variable. 
In Monte Carlo optimization, $\xi$ is usually a psuedo-random number 
generated by a computer. For most applications, such a pseudo-random 
number is sufficiently good to be regarded as a random number uniformly 
distributed on $[0, 1)$. In Section 6, we
extend the results to general situations.

To avoid repetition, we only consider the $h_{n}$ defined as in (\ref{kw_h}) 
with $\xi_{1,n}=\xi_{2,n}=\xi_{n}$. The analysis is applicable to the
one-sided finite-difference approximation (\ref{kw_1sideh}) without
any difficult. For a given $\theta_{n}$, $h_{n}$ is a finite-difference 
approximation to the derivative $J'(\theta)$ at $\theta=\theta_{n}$.
For simplicity, we omit the subscript $n$. Then
\begin{equation}
    h=\frac{L(X(\theta+\delta,\xi))-
                    L(X(\theta-\delta,\xi))}{2\delta}.
\label{h_generic}
\end{equation}
The mean of $h$ is
\[  E[h] = \frac{J(\theta+\delta)-J(\theta-\delta)}{2\delta}         \]
which is the same as that of (\ref{kw_h}) without the use of CRN. However, the
variance of (\ref{h_generic}), as we will show, is generally smaller than that of
(\ref{kw_h}) without the use of CRN when $\delta>0$ is sufficiently 
small. We will also show that the reduction in the variance 
of $h$ may have a significant impact on the convergence rate for the 
KW algorithm for Monte Carlo optimization. Toward that end, we need 
to specify the generation of the random variable $X(\theta,\xi)$ 
with a given distribution $F(\theta,x)$. Next we examine the variance 
of (\ref{h_generic}) for several popular random number generation methods. Note that
\[ Var[h]=\frac{1}{(2\delta)^{2}}
    \{E[(L(X(\theta+\delta,\xi))- L(X(\theta-\delta,\xi)))^{2}]+
         (J(\theta+\delta)- J(\theta-\delta))^{2}\}.    \]
If $J(\theta)$ is continuously differentiable on $\Theta$ with bounded 
derivatives, then 
\begin{equation}
   Var[h]=\frac{1}{(2\delta)^{2}}
    E[(L(X(\theta+\delta,\xi))- L(X(\theta-\delta,\xi)))^{2}]+O(1).
\label{h_var}
\end{equation}

\vspace{3 mm}

{\em 4.1. Inversion method}. \hspace{2mm}
Inversion is one of the most popular methods for random variable generation. 
Let $F(\theta,x)$ be the distribution function of $X(\theta,\xi)$. The
inversion method generates the random variable $X(\theta,\xi)$ in the
following way:
\begin{enumerate}
\item Generate a random number $\xi$ uniformly distributed on $[0, 1)$.
\item Set $X(\theta,\xi)=F^{-1}(\theta,\xi)$. 
\end{enumerate}
Then it is straightforward to verify that $X(\theta,\xi)$ has the desired 
distribution. Note that the mapping $F(\theta,x): R\rightarrow R$ is not one
to one in general. To ensure its existence for general 
distribution functions, the inverse function is defined as 
\[ F^{-1}(\theta,\xi)=\min\{x\;|\;F(\theta,x)> \xi,\;x\in R\} \]
which is different from the usual definition [see Krantz (1991)]. It coincides
with the usual definition if $F(\theta,x)$ is continuous and strictly
increasing. Such a definition of the inverse function covers both continuous
and discrete random variables. For example, consider a discrete random 
variable $X(\theta,\xi)=x_{i}$ with probability $p_{i}(\theta)$. 
Define $\rho_{0}(\theta)=0, \rho_{i}(\theta)=\sum_{j=1}^{i}p_{j}(\theta)$
for $i\geq 1$. Let $\xi$ be uniformly distributed on $[0, 1)$. The 
inversion method gives $F^{-1}(\theta,\xi)=x_{i}$ if $\xi\in 
[\rho_{i-1}(\theta), \rho_{i}(\theta))$.
Then direct verification shows that $X(\theta,\xi)$ obeys the desired 
distribution. This is a discrete version of the inversion method.

In order to proceed with our discussion, let us first examine the
properties of distribution functions. A distribution $F(\theta,x)$ is a 
nondecreasing and right-continuous 
function of $x$. %According to the theory of real functions,
$F(\theta,x)$ has at most countably many points of discontinuity on $R$ 
and all of the discontinuities are of the first kind --- that is, for
any $x\in R$, $F(\theta,x^{-})=\lim_{y\uparrow x}F(\theta,y)$ and 
$F(\theta,x^{+})=\lim_{y\downarrow x}F(\theta,y)$ exist and are finite 
[e.g. Krantz (1991),149-150]. Therefore, we can divide $R$ into
$\bigcup_{i}B_{i}(\theta)=R$, where $B_{i}(\theta)=[b_{i}(\theta), 
b_{i+1}(\theta))$, such that, for each $i$, $F(\theta,x)$ is continuous 
on $B_{i}(\theta)$, but jumps at $b_{i}(\theta)$. Assume that, for each 
$i$, $F(\theta,x)$ is piecewise differentiable on $B_{i}(\theta)$. Then
$F'_{x}(\theta,x)>0$ whenever it exists. We further divide the interval
$B_{i}(\theta)$ into subintervals according to whether the derivative
of $F(\theta,x)$ with respect to $x$ is zero or not. For simplicity, we 
assume that $B_{i}(\theta)= B_{i}^{0}(\theta)\bigcup B_{i}^{+}(\theta)$ 
such that $F'_{x}(\theta,x)=0$ on $B_{i}^{0}(\theta)=[b_{i}(\theta), 
c_{i}(\theta)]$ and $F(\theta,x)=F_{i}(\theta,x)$ is continuously 
differentiable with strictly positive derivatives on $B_{i}^{+}(\theta)=
(c_{i}(\theta), b_{i+1}(\theta))$. It is possible that $b_{i}(\theta)=
c_{i}(\theta)$. On $B_{i}^{0}(\theta)$, the derivatives $F'_{x}(\theta,x)$ 
should be understood as the right and the left derivatives at 
$b_{i}(\theta), c_{i}(\theta)$, respectively. It is possible that
$F(\theta,x)$ is not differentiable at $c_{i}(\theta)$. The inverse 
$F_{i}^{-1}(\theta,\xi)$ is defined in the usual sense. It is continuous, 
strictly increasing, and differentiable on $(F(\theta,c_{i}(\theta)), 
F(\theta,b_{i+1}^{-}(\theta)))$. 

Under the preceding decomposition, $F(\theta,x)$ is
discontinuous at $b_{i}(\theta)$, is a constant on $B_{i}^{0}(\theta)$, 
and is strictly increasing and differentiable on $B_{i}^{+}(\theta)$.

The following Lemma 2 follows directly from the definition of the inverse
function and the decomposition of $F(\theta,x)$.
\vspace{3 mm}

{\sc Lemma 2}. {\em Let $X(\theta,\xi)$ be defined by the inverse function
$X(\theta,\xi)=F^{-1}(\theta,\xi)$. Let $\Xi_{i}(\theta)=
[F(\theta,b_{i}^{-}(\theta)), F(\theta,b_{i+1}^{-}(\theta)))$. 
Then $\Xi_{i}(\theta)\subset[0, 1)$ and for any $\xi\in \Xi_{i}(\theta)$
\begin{equation}
  X(\theta,\xi)=\left\{ \begin{array}{ll}
      b_{i}(\theta), & \mbox{if } \xi\in [F(\theta,b_{i}^{-}(\theta)),
                                 F(\theta,c_{i}(\theta))), \\
      c_{i}(\theta), & \mbox{if } \xi=F(\theta,c_{i}(\theta)), \\
      F_{i}^{-1}(\theta,\xi), & \mbox{if } \xi\in (F(\theta,c_{i}(\theta)),
                                 F(\theta,b_{i+1}^{-}(\theta))).
     \end{array} \right.
\label{X1}
\end{equation}}

\vspace{4 mm}
We need the following result.
\vspace{3 mm}

{\sc Lemma 3}. {\em Assume that
\begin{itemize}
\item[(C1).] $L(X)$ and $L'_{X}(X)$ are bounded, $J(\theta)$ is continuously 
             differentiable on $\Theta$;
\item[(C2).] for each $i$, $F_{i}(\theta,x)$ is continuously differentiable
            on $B_{i}^{+}(\theta)$ with strictly positive derivatives with
            respect to $x$, and 
       \[  \sum_{i}E[(\max_{\theta}(F'_{i\theta}(\theta,x))^{2}/
                 F'_{i x}(\theta,x))I_{B_{i}^{+}(\theta)}]<\infty;  \]
\item[(C3).] $b_{i}(\theta)$ is continuously differentiable in $\theta$, and
      \( \sum_{i}\max_{\theta}(b_{i}'(\theta))^{2}< \infty; \)
\item [(C4).] for each $i$, the functions $F(\theta,c_{i}(\theta))$ and 
             $F(\theta,b_{i}^{-}(\theta))$ are continuously 
             differentiable in $\theta$, and
      \( \sum_{i}\max_{\theta}|F'(\theta,c_{i}(\theta))|<\infty,\)
   \( \sum_{i}\max_{\theta}|F'(\theta,b_{i}^{-}(\theta))|<\infty,\)
\end{itemize}
Define
\(   M_{1}(\theta) = 2\sum_{i}(L(c_{i}(\theta))-L(b_{i}(\theta)))^{2}
           |F'(\theta,c_{i}(\theta))|.       \)
Then $M_{1}(\theta)\geq 0$ is bounded for all $\theta$. 
If $M_{1}(\theta)>0$, we have
\begin{equation}
      E[(L(X(\theta+\delta,\xi))- L(X(\theta-\delta,\xi)))^{2}]
       =M_{1}(\theta)\delta+o(\delta)
\label{EFD}
\end{equation}
as $\delta>0$ goes to zero}.

{\sc Proof}. We calculate 
\begin{eqnarray*}
 \lefteqn{\lim_{\delta\rightarrow 0}
         \frac{1}{\delta}E[(L(X(\theta+\delta,\xi))
         -L(X(\theta-\delta,\xi)))^{2}]} \\
  & & = \lim_{\delta\rightarrow 0}
       \frac{1}{\delta}\sum_{i}E[(L(X(\theta+\delta,\xi))
      -L(X(\theta-\delta,\xi)))^{2}I_{\Xi_{i}(\theta-\delta)}].
\end{eqnarray*}
Let
\[  R_{i}(\theta,\delta)=\frac{1}{\delta}E[(L(X(\theta+\delta,\xi))
      -L(X(\theta-\delta,\xi)))^{2}I_{\Xi_{i}(\theta-\delta)}].  \]
Then 
\begin{equation}
  \lim_{\delta\rightarrow 0}\frac{1}{\delta}E[(L(X(\theta+\delta,\xi))
    -L(X(\theta-\delta,\xi)))^{2}] = \lim_{\delta\rightarrow 0}
       \sum_{i} R_{i}(\theta,\delta).
\label{fd1}
\end{equation}
Next, we prove that the limit and the summation commute.
Define $D_{i,1} = \Xi_{i}(\theta-\delta)\bigcap [0, F(\theta+\delta,
b_{i}^{-}(\theta+\delta)))$, $D_{i,2} = \Xi_{i}(\theta-\delta)\bigcap 
\Xi_{i}(\theta+\delta)$, and $D_{i,3} = \Xi_{i}(\theta-\delta)\bigcap 
[F(\theta+\delta,b_{i+1}^{-}(\theta+\delta)), 1)$. It is possible for each of 
$D_{i,j}, j=1, 2,3,$ to be empty. Then $\Xi_{i}(\theta-\delta)=\bigcup
D_{i,j}$ and
\begin{equation}
  R_{i}(\theta,\delta) = \sum_{j=1}^{3} R_{i,j}, \;\; R_{i,j}=
      \frac{1}{\delta}E[(L(X(\theta+\delta,\xi))
      -L(X(\theta-\delta,\xi)))^{2}I_{D_{i,j}}].   
\label{r1}
\end{equation}
By Assumption (C1), there exist $N_{1}, N_{2}>0$ such that $|L(X)|\leq 
N_{1}$, $|L'_{X}(X)|\leq N_{2}$. Therefore,
\begin{eqnarray*}
    R_{i,1} & \leq & (2N_{1})^{2}|F(\theta+\delta,b_{i}^{-}(\theta+\delta))-
          F(\theta-\delta,b_{i}^{-}(\theta-\delta))|/\delta      \\
     & \leq & 2(2N_{1})^{2}\max_{\theta}|F'(\theta,
             b_{i}^{-}(\theta))|,                       
\end{eqnarray*}
\begin{eqnarray*}
   R_{i,3} & \leq & (2N_{1})^{2}|F(\theta+\delta,b_{i+1}^{-}(\theta+\delta))-
          F(\theta-\delta,b_{i+1}^{-}(\theta-\delta))|/\delta     \\
      & \leq & 2(2N_{1})^{2}\max_{\theta}|F'(\theta,
             b_{i+1}^{-}(\theta))|,               
\end{eqnarray*}
Without loss of generality, assume that $F(\theta+\delta,b_{i}^{-}(\theta+
\delta)))\geq F(\theta-\delta,b_{i}^{-}(\theta-\delta)))$ and $F(\theta+\delta,
b_{i+1}^{-}(\theta+\delta)))\leq F(\theta-\delta,b_{i+1}^{-}(\theta-\delta)))$.
If $F(\theta+\delta,c_{i}(\theta+\delta))
>F(\theta-\delta,c_{i}(\theta-\delta))$, 
\begin{eqnarray}
  R_{i,2} & = & \frac{1}{\delta}\int_{F(\theta+\delta,b_{i}^{-}(\theta
            +\delta))}^{F(\theta-\delta,c_{i}(\theta-\delta))} 
            (L(X(\theta+\delta,\xi))
      -L(X(\theta-\delta,\xi)))^{2}d\xi   
\label{rr2}
\end{eqnarray}
\begin{eqnarray*}         
  & & \mbox{} + \frac{1}{\delta}\int_{F(\theta-\delta,c_{i}(\theta-\delta))}
         ^{F(\theta+\delta,c_{i}(\theta+\delta))} 
    (L(X(\theta+\delta,\xi))-L(X(\theta-\delta,\xi)))^{2}d\xi \nonumber \\
 & & \mbox{} + \frac{1}{\delta}\int_{F(\theta+\delta,c_{i}(\theta+\delta))}
        ^{F(\theta+\delta,b_{i+1}^{-}(\theta+\delta))} 
    (L(X(\theta+\delta,\xi))-L(X(\theta-\delta,\xi)))^{2}d\xi \nonumber \\
 & = & \frac{1}{\delta}\int_{F(\theta+\delta,b_{i}^{-}(\theta
         +\delta)))}^{F(\theta-\delta,c_{i}(\theta-\delta))} 
      (L(b_{i}(\theta+\delta))-L(b_{i}(\theta-\delta)))^{2}d\xi  \nonumber \\ 
  & & \mbox{} + \frac{1}{\delta}\int_{F(\theta-\delta,c_{i}(\theta-\delta))}
         ^{F(\theta+\delta,c_{i}(\theta+\delta))} (L(b_{i}(\theta+\delta))
         -L(F_{i}^{-1}(\theta-\delta,\xi)))^{2}d\xi    \nonumber \\
 & & \mbox{} + \frac{1}{\delta}\int_{F(\theta+\delta,c_{i}(\theta+\delta))}
          ^{F(\theta+\delta,b_{i+1}^{-}(\theta+\delta))} 
        (L(F_{i}^{-1}(\theta+\delta,\xi))
         -L(F_{i}^{-1}(\theta-\delta,\xi)))^{2}d\xi  % \nonumber
%\label{rr2}
\end{eqnarray*}
The first two terms of (\ref{rr2}) are bounded respectively by 
\[   4N_{2}^{2}\max_{\theta}(b'_{i}(\theta))^{2}\delta\;\;\mbox{ and }\;\;
       2(2N_{1})^{2}\max_{\theta}|F'(\theta,c_{i}(\theta))|.  \] 
The third term of (\ref{rr2}) can be rewritten as
\[  \frac{1}{\delta}\int_{F(\theta+\delta,c_{i}(\theta+\delta))}
          ^{F(\theta,c_{i}(\theta))} 
        (L(F_{i}^{-1}(\theta+\delta,\xi))
         -L(F_{i}^{-1}(\theta-\delta,\xi)))^{2}d\xi    \]
\[   + \frac{1}{\delta}\int_{F(\theta,b_{i+1}^{-}(\theta))}
          ^{F(\theta+\delta,b_{i+1}^{-}(\theta+\delta))} 
        (L(F_{i}^{-1}(\theta+\delta,\xi))
         -L(F_{i}^{-1}(\theta-\delta,\xi)))^{2}d\xi        \]
\[  + \frac{1}{\delta}\int_{F(\theta,c_{i}(\theta))}
          ^{F(\theta,b_{i+1}^{-}(\theta))} 
        (L(F_{i}^{-1}(\theta+\delta,\xi))
         -L(F_{i}^{-1}(\theta-\delta,\xi)))^{2}d\xi        \]
\[  \leq 2(2N_{1})^{2}\max_{\theta}|F'(\theta,c_{i}(\theta))|
    +2(2N_{1})^{2}\max_{\theta}|F'(\theta,b_{i+1}^{-}(\theta))|  \]
\[ + 4N_{2}^{2}E[(\max_{\theta}(F'_{i\theta}(\theta,x))^{2}/
          F'_{i x}(\theta,x))I_{B_{i}^{+}(\theta)}]\delta.  \]
Therefore, $R_{i,2}$ is bounded by
\[   4N_{2}^{2}\max_{\theta}(b'_{i}(\theta))^{2}\delta+
    4(2N_{1})^{2}\max_{\theta}|F'(\theta,c_{i}(\theta))|  \]
\[  +2(2N_{1})^{2}\max_{\theta}|F'(\theta,b_{i+1}^{-}(\theta))|  
  + 4N_{2}^{2}E[(\max_{\theta}(F'_{i\theta}(\theta,x))^{2}/
      F'_{i x}(\theta,x))I_{B_{i}^{+}(\theta)}]\delta  \]
Similarly, we can prove that if $F(\theta+\delta,c_{i}(\theta+\delta))
\leq F(\theta-\delta,c_{i}(\theta-\delta))$,
\begin{eqnarray}
  R_{i,2} & = & \frac{1}{\delta}\int_{F(\theta+\delta,b_{i}^{-}(\theta
            +\delta)))}^{F(\theta+\delta,c_{i}(\theta+\delta)} 
      (L(b_{i}(\theta+\delta))-L(b_{i}(\theta-\delta)))^{2}d\xi 
 \label{rr3}
 \end{eqnarray}
 \begin{eqnarray*}
  & & \mbox{} + \frac{1}{\delta}\int_{F(\theta+\delta,c_{i}(\theta+\delta))}
         ^{F(\theta-\delta,c_{i}(\theta-\delta))}
         (L(F_{i}^{-1}(\theta+\delta,\xi))
        -L(b_{i}(\theta-\delta)))^{2}d\xi \nonumber \\
 & & \mbox{} + \frac{1}{\delta}\int_{F(\theta-\delta,c_{i}(\theta-\delta))}
        ^{F(\theta+\delta,b_{i+1}^{-}(\theta+\delta)))} 
            (L(F_{i}^{-1}(\theta+\delta,\xi))
           -L(F_{i}^{-1}(\theta-\delta,\xi)))^{2}d\xi \nonumber \\
 &\leq & 4N_{2}^{2}\max_{\theta}(b'_{i}(\theta))^{2}\delta 
  +4(2N_{1})^{2}\max_{\theta}|F'(\theta,c_{i}(\theta))|  \\
  & &\mbox{}+2(2N_{1})^{2}\max_{\theta}
              |F'(\theta,b_{i+1}^{-}(\theta))|\nonumber \\
  & & \mbox{} + 4N_{2}^{2}E[(\max_{\theta}(F'_{i\theta}(\theta,x))^{2}/
      F'_{ix}(\theta,x))I_{B_{i}^{+}(\theta)}]\delta. \nonumber
%\label{rr3}
\end{eqnarray*}
Substituting the upper bounds for $R_{i,j}, j=1,2,3,$ into
(\ref{fd1}), also noting the assumptions (C2)-(C4), we see that 
$R_{i}(\theta,\delta)$ is uniformly bounded with respect to $\delta$. 
Therefore, $\sum_{i}R_{i}(\theta,\delta)$ converges uniformly
in $(0,\delta_{0})$ for any $\delta_{0}>0$. By the Weierstrass M-test 
[Krantz (1991),211], we know that
the limit and the summation commute. From (\ref{fd1}), 
\begin{equation}
  \lim_{\delta\rightarrow 0}\frac{1}{\delta}E[(L(X(\theta+\delta,\xi))
    -L(X(\theta-\delta,\xi)))^{2}] = \lim_{\delta\rightarrow 0}
       \sum_{i} R_{i}(\theta,\delta) 
     = \sum_{i} \lim_{\delta\rightarrow 0}R_{i}(\theta,\delta).
\label{fd2}
\end{equation}
We next calculate $\lim_{\delta\rightarrow 0}
R_{i}(\theta,\delta)$. For each $i$, there exists a $\delta_{i}>0$ 
such that for any $\delta\leq \delta_{i}$
\begin{eqnarray*}
  \lefteqn{ |F(\theta+\delta,b_{j}^{-}(\theta+\delta)))- 
    F(\theta-\delta,b_{j}^{-}(\theta-\delta)))|} \\
 & & \leq \frac{1}{4}\min_{j=i-1,i,i+1,i+2}
       \{F(\theta-\delta,b_{j+1}^{-}(\theta-\delta))
       -F(\theta-\delta,b_{j}^{-}(\theta-\delta))\}. 
\end{eqnarray*}
Note that $D_{i,1}=\Xi_{i}(\theta-
\delta)\bigcap\Xi_{i-1}(\theta+\delta)$ and $D_{i,3}=\Xi_{i}(\theta-
\delta)\bigcap\Xi_{i+1}(\theta+\delta)$ when $\delta\leq\delta_{i}$.
Therefore, by taking into account that each of $D_{i,1}$ and $D_{i,3}$
may be empty, we have
\begin{eqnarray*}
    R_{i,1} & \leq & \int^{F(\theta+\delta,b_{i}^{-}(\theta+\delta))}
          _{F(\theta-\delta,b_{i}^{-}(\theta-\delta))}     
            (L(b_{i}(\theta+\delta))-L(F_{i-1}^{-1}(\theta
            -\delta,\xi)))^{2}d\xi               \\
   & \leq & \max_{\theta}|F'(\theta,b_{i}^{-}(\theta))|
  (L(b_{i}(\theta+\delta))-L(F_{i-1}^{-1}(\theta-\delta,\tilde{\xi})))^{2} \\
    & = & o(1) 
\end{eqnarray*}
where $\tilde{\xi}\in[F(\theta-\delta,b_{i}^{-}(\theta-\delta)),
F(\theta+\delta,b_{i}^{-}(\theta+\delta)))$. Similarly, $R_{i,3} = o(1)$.
Hence, $\lim_{\delta \rightarrow 0} D_{i,j} = 0$ for $j=1, 3$. 
Also, the analysis of (\ref{rr2}) and (\ref{rr3}) shows that
\begin{equation}
   \lim_{\delta\rightarrow 0} D_{i,2}
  = 2(L(c_{i}(\theta))-L(b_{i}(\theta)))^{2}
         |F'(\theta,c_{i}(\theta))|
\label{fd3}
\end{equation}
Substituting (\ref{fd3}) into (\ref{fd2}) we get
\[   \lim_{\delta\rightarrow 0}\frac{1}{\delta}E[(L(X(\theta+\delta,\xi))
    -L(X(\theta-\delta,\xi)))^{2}] = M_{1}(\theta)  \]
which is exactly what we want to prove.
\hspace{0.1in}\rule{1.5mm}{3.0mm}
\vspace{3 mm}

The proof of Lemma 3 shows that Assumptions (C2) and (C3)
guarantee that the inverse function $F_{i}^{-1}(\theta,\xi)$ is
sufficiently smooth. Assumption (C4) ensures the existence of 
$M_{1}(\theta)$. Assumptions (C2)-(C4) are mild. Assumption (C1)
guarantees the smoothness of the function $L(X)$. The boundedness of
$L(X)$ and $L'_{X}(X)$ can be removed if there are only a finite number 
of sets of $B_{i}(\theta)$. The finiteness of $B_{i}(\theta)$ can also 
relax the assumptions (C2)-(C4).

The case of $M_{1}(\theta) = 0$ can only occur when either $b_{i}(\theta)=
c_{i}(\theta)$ or $F'(\theta,c_{i}(\theta)) =0$. The situation 
of $b_{i}(\theta)=c_{i}(\theta)$ (assuming that $L(X)$ is not a constant) 
happens when $F(\theta,x)$ is strictly increasing. A repetition of the 
proof of Lemma 3 yields that
\vspace{3 mm}

{\sc Corollary 5}. {\em If Assumption (C1) is satisfied and
\begin{equation}
    E[(F'_{\theta}(\theta,x))^{2}/F'_{x}(\theta,x)]<\infty,   
\label{clr5_1}
\end{equation}
then $E[(L(X(\theta+\delta,\xi))-L(X(\theta-\delta,\xi)))^{2}]=
O(\delta^{2})$.}
\vspace{3 mm}

Corollary 5 recovers a result obtained by Glasserman and Yao (1992) 
under the assumption of Lipschitz continuity of $L(F^{-1}(\theta,\xi))$.
When $F'(\theta,c_{i}(\theta)) =0$ for all $i$, using the same
arguments as that of Corollary 5 we can establish that
\vspace{3 mm}

{\sc Corollary 6}. {\em In addition to Assumptions (C1)-(C4), assume 
that $F(\theta,c_{i}(\theta))$ is continuously twice differentiable for 
all $i$ with 
\begin{equation}
   0< \sum_{i}(L(c_{i}(\theta))-L(b_{i}(\theta)))^{2}
    |F''(\theta,c_{i}(\theta))|<\infty.   
\label{clr6_1}
\end{equation}
Then, $E[(L(X(\theta+\delta,\xi))-L(X(\theta-\delta,\xi)))^{2}]=
O(\delta^{2})$.}

\vspace{3 mm}

The following Theorem 4 is the main conclusion of this subsection.
\vspace{3 mm}

{\sc Theorem 4}. {\em Assume that Assumptions (A1)-(A4) and (C1)-(C4)
are satisfied. If $M_{1}(\theta) > 0$ for all $\theta$, 
then the best convergence rate 
for the KW algorithm (\ref{kw}) with $h_{n}$ defined by (\ref{h_generic}) is $n^{-2/5}$. 
This rate is attained by choosing $a_{n}=an^{-1}$, $a>2/(5K_{1})$, and 
$\delta_{n}=n^{-1/5}$.} 

{\sc Proof}. Under Assumption of (C1)-(C4) and $M_{1}(\theta)>0$, we 
know from Lemma 3 that $E[(L(X(\theta+\delta,\xi))-L(X(\theta-\delta,\xi)))^{2}]
=M_{1}(\theta)\delta+o(\delta)$. According to (\ref{h_var}), the variance of $h_{n}$
is of order $Var[h_{n}|\theta_{n}]=M_{1}(\theta_{n})/\delta_{n}
+o(1/\delta_{n})$. Lemma 3 shows that $M_{1}(\theta)$ is bounded.
Therefore, $\gamma=-1$ in (\ref{h2}). Since (A1)-(A4) are satisfied, Corollary 2 
shows that the best convergence rate is $n^{-2/(4-\gamma)}=n^{-2/5}$.
\hspace{0.1in}\rule{1.5mm}{3.0mm}
\vspace{3 mm}

The following Theorem 5 summerizes the rate of convergence
of the KW algorithm (\ref{kw}) when (\ref{kw_h}) is replaced with
one-sided finite difference approximation with CRN. 
The proofs are omitted since they are very
similar to that of Theorem 4.
\vspace{3 mm}

{\sc Theorem 5}. {\em (I) Under the same conditions as those of Theorem 4 but
the estimate $h$ is replaced by the following one-sided finite difference
with the use of CRN
\begin{equation}
    h=\frac{L(X(\theta+\delta,\xi))- L(X(\theta,\xi))}{\delta},
\label{th5_1}
\end{equation}
the best convergence rate is $n^{-1/3}$ which is achieved by setting
$a_{n}=an^{-1}, a>1/3K_{1}$, and $\delta_{n}=d n^{-1/3}$.

(II) Assume all the assumptions of Theorem 4 except that $M_{1}(\theta) 
= 0$. Then Corollaries 5 and 6 show that $E[(L(X(\theta+\delta,\xi))
-L(X(\theta-\delta,\xi)))^{2}]= O(\delta^{2})$ if either of (\ref{clr5_1}) or 
(\ref{clr6_1}) holds. Hence, 
$Var[h_{n}|\theta_{n}]=O(1)$ for $h_{n}$ defined 
by either (\ref{h_generic}) or (\ref{th5_1}). The best convergence rate for the 
KW algorithm (\ref{kw})
is $n^{-1/2}$. This rate can be attained by setting $a_{n}=an^{-1}, 
a>1/2K_{1}$, and $\delta_{n} = d n^{-\eta}, \eta\geq 1/2$.}
\vspace{3 mm}

We would like to emphasize that the assumptions in Corollaries 5 and 6
are satisfied for a broad class of stochastic optimization problems
[see Glasserman and Yao (1992) for a discussion].
Theorems 4 and 5 state that, when the inversion method is used
in the generation of random variables and when $h$ is defined by (\ref{h_generic}), 
the convergence rate for the KW algorithm with CRN is $n^{-2/5}$ in 
general and is $n^{-1/2}$ for a large class of problems that satisfy the 
assumptions in Corollaries 5 and 6. The improvement is signficant since the
best possible rate for the same KW algorithm without CRN is $n^{-1/3}$.
\vspace{3 mm}

{\em 4.2. Rejection method}.\hspace{2mm}
Let $f(\theta,x)$ be the density function of $X(\theta,\xi)$. 
Assume that, for all $\theta\in \Theta$, $f(\theta,x)$ is zero outside 
a finite interval $[a, b]$ and is bounded by $0\leq f(\theta,x)
\leq c$, $c>0$ is a constant. The rejection method
generates $X(\theta,\xi)$ according to the following three steps:
\begin{enumerate}
\item Generate $\xi_{1}$ uniformly distributed on $[a, b]$.
\item Generate $\xi_{2}$ uniformly distributed on $[0, c]$.
\item If $\xi_{2}\leq f(\theta,\xi_{1})$, then set $X(\theta,\xi)=\xi_{1}$;
      otherwise go to 1.
\end{enumerate}
The rejection method uses at least two random numbers $\xi_{1}$
and $\xi_{2}$ to generate $X(\theta,\xi)$. The total number of random 
numbers $\xi_{1}, \xi_{2}$ required before outputing $X(\theta,\xi)$ 
is a random value. The rejection method does not accurately meet the 
CRN requirements since it is impossible to define $X(\theta+\delta,\xi)$ 
and $X(\theta-\delta,\xi)$ using a fixed set of
uniform random numbers [Bratley et al. (1983); Franta (1975)]. 
Therefore, we modify the definition of CRN in the sense defined by the
following procedure for the generation of a paired random variables: 
\vspace{3 mm}

\noindent {\em Generation of $X(\theta+\delta,\xi)$ and 
$X(\theta-\delta,\xi)$}:
\begin{enumerate}
\item Generate $\xi_{1}$ uniformly distributed on $[a, b]$.
\item Generate $\xi_{2}$ uniformly distributed on $[0, c]$.
\item If $\xi_{2}\leq f(\theta-\delta,\xi_{1})$ and $\xi_{2}\leq 
      f(\theta+\delta,\xi_{1})$, then set $X(\theta-\delta,\xi)=
      X(\theta+\delta,\xi) =\xi_{1}$.
\item If $\xi_{2}\leq f(\theta-\delta,\xi_{1})$ and $\xi_{2}>
      f(\theta+\delta,\xi_{1})$, then set $X(\theta-\delta,\xi) 
      =\xi_{1}$ and generate a $X(\theta+\delta,\xi)=\xi_{3}$ by the 
      rejection method.
\item If $\xi_{2} > f(\theta-\delta,\xi_{1})$ and $\xi_{2}  \leq
      f(\theta+\delta,\xi_{1})$, then set $X(\theta+\delta,\xi) 
      =\xi_{1}$ and generate a $X(\theta-\delta,\xi)=\xi_{4}$ by the 
      rejection method.
\item If $\xi_{2}> f(\theta-\delta,\xi_{1})$ and $\xi_{2}>
      f(\theta+\delta,\xi_{1})$, go to 1.
\end{enumerate}
This is essentially a coupling procedure [see Devroye (1990) for a 
discussion on coupling]. 
Such a modification is necessary to mimic the scheme of CRN using the 
rejection method. We will soon see that even such a loosely defined 
scheme can accelerate the convergence of the KW algorithm.
Let $X(\theta-\delta,\xi), X(\theta+\delta,\xi)$ be generated by the 
preceding procedure. It is obvious that $E[h]$ for $h$ in (\ref{h_generic}) remains 
the same as that in the inversion method. 
\vspace{3 mm}

{\sc Theorem 6}. {\em Suppose that $f(\theta,x)$ is zero outside $[a, b]$, 
$0\leq f(\theta,x)\leq c$ for all $x\in [a, b]$, and $X(\theta-\delta,\xi),
X(\theta+\delta,\xi)$ are generated by the previously described 
procedure. Assume that 
\begin{itemize}
\item[(H1).] $Var[L(X(\theta,\xi))]$ is continuous in $\theta\in\Theta$;
\item[(H2).] $f(\theta,x)$ is differentiable in $\theta$ for each 
             $x\in [a, b]$, $f(\theta,x)$ satisfies the Lipschitz 
             condition with respect to $\theta$, i.e., there is a
             $K(x)$ such that $|f(\theta+\delta,x)-f(\theta,x)|\leq
             K(x)\delta$, and that $\int_{a}^{b}K(x)dx<\infty$. 
\end{itemize}
Define
\[   M_{2}(\theta) = \frac{Var[L(X(\theta,\xi))]}{2c(b-a)}
     \int_{a}^{b}|f'_{\theta}(\theta,x)|dx.   \]
Then $0\leq M_{2}(\theta)<\infty$. If $M_{2}(\theta)>0$ for all $\theta$,
$Var[L(X(\theta,\xi))]$ is bounded,
$h$ is defined by (\ref{h_generic}), and the assumptions (A1)-(A4) are satisfied, then 
the convergence rate for the KW algorithm with CRN is $n^{-2/5}$}.

{\sc Proof}. We see from the procedure of generating $X(\theta-\delta,\xi)$
and $X(\theta+\delta,\xi)$ that, conditioned on either $\xi_{2}\leq
f(\theta-\delta,\xi_{1})$ or $\xi_{2}\leq f(\theta+\delta,\xi_{1})$,
$X(\theta+\delta,\xi)= X(\theta-\delta,\xi)=\xi_{1}$ when $\xi_{2}\leq
f(\theta-\delta,\xi_{1})$ and $\xi_{2}\leq f(\theta+\delta,\xi_{1})$;
otherwise $X(\theta+\delta,\xi)=\xi_{3}$ and $X(\theta-\delta,\xi)=\xi_{4}$. 
Note that $\xi_{3}$ and $\xi_{4}$ are independent. Therefore,
\begin{equation}
    Var[h] = \frac{1}{4\delta^{2}}(Var[L(\xi_{3})]
         +Var[L(\xi_{4})])
       \frac{1}{c}E[|f(\theta+\delta,\xi_{1})-f(\theta-\delta,\xi_{1})|]
\label{th6_varh}
\end{equation}
Under Assumption (H1), $Var[L(\xi_{3})]+Var[L(\xi_{4})]= 
2Var[L(X(\theta,\xi))]+o(1)$. By Assumption (H2), 
\[ \frac{1}{\delta}E[|f(\theta+\delta,\xi_{1})-f(\theta-\delta,\xi_{1})|]
     \leq\frac{2}{b-a}\int_{a}^{b}K(x)dx   \]
and $K(x)$ is integrable on $[a, b]$. According to the Weierstrass M-test,
(\ref{th6_varh}) implies that 
\begin{eqnarray*}
    Var[h] &=& \frac{1}{2\delta^{2}}Var[L(X(\theta,\xi))]
       \frac{1}{c}E[|f(\theta+\delta,\xi_{1})-f(\theta-\delta,\xi_{1})|]
                   +o(\frac{1}{\delta^{2}})          \nonumber \\
    &  = & \frac{1}{2\delta}Var[L(X(\theta,\xi))]
       \frac{1}{c(b-a)}\int_{a}^{b}|
   f'_{\theta}(\theta,x)|dx+o(\frac{1}{\delta}) \nonumber \\    
    &  = &\frac{M_{2}(\theta)}{\delta}+o(\frac{1}{\delta}).   
\end{eqnarray*}
Thus, we know from Corollary 2 where $\gamma=-1$ that the 
conclusion follows.
\hspace{0.1in}\rule{1.5mm}{3.0mm}
\vspace{3 mm}

For simplicity, we only consider the simplest form of the rejection method
and the case in which $f(\theta,x)$ is continuous. An analysis similar
to the one used in the proof of Theorem 6 shows that 
$Var[h]=O(1/\delta)$ remains valid in the following three situations:
(i) The estimate $h$ is replaced by the one-sided finite difference (\ref{th5_1}); 
(ii) The density function $f(\theta,x)$ is piecewise differentiable; (iii) The
rejection method is replaced by the following {\em generalized rejection
method}. Assume that there exist a positive constant $A$ and a density
function $g(x)$ such that $f(\theta,x)\leq A g(x)$ for all $\theta$
and for all $x\in [a, b]$. Then
\begin{enumerate}
\item generate $\xi_{1}$ with the density function $g(x)$;
\item generate $\xi_{2}$ uniformly distributed on $[0, A g(\xi_{1})]$;
\item if $\xi_{2}\leq f(\theta,\xi_{1})$, then set $X(\theta,\xi)=\xi_{1}$;
      otherwise go to 1.
\end{enumerate}
It is easy to verify that $X(\theta,\xi)$ has the desired distribution.
The density function $g(x)$ should be chosen such that it is easier to 
generate a random variable with $g(x)$ than those with $f(\theta,x)$.

Generally speaking, the convergence rates for the KW algorithm are the
same when either the inversion method or the rejection method is used
in the generation of the random variable $X(\theta,\xi)$. However, the 
rate corresponding to the use of the rejection method is universally
true for any function: It can be seen from its definition that
$M_{2}(\theta)$ is always positive except when $Var[L(X)]=0$ 
or when $f(\theta,x)$ is
independent of $\theta$. Both cases are of little practical relevance. 
Furthermore, assume that the assumptions in Theorem 6 are 
satisfied and, in addition, $f(\theta,x)$ is strictly positive on
$(a, b)$. Then the best possible convergence rate for the KW algorithm 
is $n^{-2/5}$ when the rejection method is used in generating 
$X(\theta,\xi)$. On the other hand, the distribution 
function $F(\theta,x)=\int_{a}^{x}f(\theta,u)du$ is continuously 
differentiable and strictly increasing on $[a, b]$. Theorem 5
shows that the convergence rate for the KW algorithm is $n^{-1/2}$ if the
inversion method is used in generating $X(\theta,\xi)$. Therefore, as far
as the convergence of the KW algorithm is concerned, the inversion method
leads to faster convergence than the rejection method. This conclusion is in favor of
the argument that the inversion method is superior to the rejection 
method [c.f. Bratley et al. (1983), 141].
\vspace{3 mm}

{\em 4.3. Composition method}.\hspace{2mm}
Assume that the distribution function $F(\theta,x)$ of $X(\theta,\xi)$ is
of the form
\[    F(\theta,x) = \sum_{i=1}^{m} p_{i}(\theta)F_{i}(\theta,x)  \]
where $p_{i}(\theta)>0, m\leq\infty, \sum_{i}p_{i}(\theta)=1$, and for
each $i$, $F_{i}(\theta,x)$ is a distribution function. The composition 
method generates the random variable $X(\theta,\xi)$ in the following way:
\begin{enumerate}
\item Generate a random variable $Y$ with distribution $Prob\{Y=i\}
      =p_{i}(\theta)$.
\item If $Y=i$, generate $X(\theta,\xi)$ according to distribution 
      $F_{i}(\theta,x)$.
\end{enumerate}

In the composition method, there is no specification on the method for
the generation
of random variables at each step. Any method such as inversion and 
rejection can be used. As an example, we consider the case in
which random variables are generated using the inversion method which
is superior to the rejection method, as we have argued in the 
previous subsection. Define $\rho_{0}(\theta)=0, \rho_{i}(\theta)
=\sum_{j=1}^{i}p_{j}(\theta)$ for $i\geq 1$. The following procedure 
is the actual composition method we are considering.
\begin{enumerate}
\item Generate a random number $\xi_{1}$ uniformly distributed on $[0, 1)$.
\item If $\xi_{1}\in[\rho_{i-1}(\theta),\rho_{i}(\theta))$, then generate
      a random number $\xi_{2}$ uniform on $[0, 1)$ and set $X(\theta,\xi)=
      F_{i}^{-1}(\theta,\xi_{2})$.
\end{enumerate}

In this algorithm, we need two uniform random numbers in the generation
of $X(\theta,\xi)$. Actually we can do with only one random number by setting
$\xi_{2}=(\xi_{1}-\rho_{i-1}(\theta))/p_{i}(\theta)$. Direct verification shows
that, conditional on $\xi_{1}\in [\rho_{i-1}(\theta),\rho_{i}(\theta))$, 
$\xi_{2}$ is uniform on $[0, 1)$. In the composition
method, we regard that $X(\theta-\delta,\xi)$ and $X(\theta+\delta,\xi)$
conform the CRN requirement if they are generated by the preceding 
procedure using the same $\xi=(\xi_{1},\xi_{2})$. We can prove that 
it can accelerate the convergence of the KW algorithm.

For simplicity, we assume that, for each $i$, $F_{i}(\theta,x)
=F_{i}(x)$ is independent of $\theta$, and  the number of distribution
component is finite, i.e., $m<\infty$. The case in which
$F(\theta,x)$ is of general form can be treated parallel to the proof
of Theorem 4. Our aim here is to find special features of the
decomposition method rather than to develop the complete theory 
which is not difficult
to derive. We first consider the situation where $\xi_{1}$ and 
$\xi_{2}$ are independent. We then examine the case where $\xi_{2}=
(\xi_{1}-\rho_{i-1}(\theta))/p_{i}(\theta)$.
\vspace{3 mm}

{\sc Theorem 7}. {\em Suppose that $\xi_{1}$ and $\xi_{2}$ are independent
, $p_{i}(\theta)$ is differentiable in $\theta$, and
\[      M_{3}(\theta) = \sum_{i=1}^{m}E[(L(F_{i+1}^{-1}(\xi))
               -L(F_{i}^{-1}(\xi)))^{2}]|\rho'_{i}(\theta)|        \]
exits and is finite. If $M_{3}(\theta)>0$ is bounded from above
for all $\theta$, $h$ is defined
by (\ref{h_generic}), and the assumptions (A1)-(A4) are satisfied, then the convergence 
rate for the KW algorithm is $n^{-2/5}$.}

{\sc Proof}. The proof is a simplified version of that of Lemma 3 since the 
assumptions here are stronger. According to the generation scheme of $X(\theta+\delta,\xi)$ and $X(\theta-\delta,\xi)$, we have 
\begin{eqnarray}
   E[(L(X(\theta+\delta,\xi))
               -L(X(\theta-\delta,\xi)))^{2}]  
\label{th7_1}
\end{eqnarray}
\begin{eqnarray*}               
 & & =\sum_{i=1}^{m}\int_{\rho_{i-1}(\theta-\delta)}^{\rho_{i}(\theta-\delta)}
   E_{\xi_{2}}[(L(X(\theta+\delta,\xi))-
       L(X(\theta-\delta,\xi)))^{2}]d\xi_{1}      \nonumber
\end{eqnarray*}
Note that $m$ is finite and, for each $i$, $\rho_{i}(\theta)$ is continuous.
There exists a $\delta_{0}$ such that when $\delta\leq\delta_{0}$
\begin{equation}
     |\rho_{i}(\theta+\delta)-\rho_{i}(\theta-\delta)|
       \leq \frac{1}{4}\min_{j}\{\rho_{j+1}(\theta-\delta)
                -\rho_{j}(\theta-\delta)\}, \;\;\mbox{ for all $i$.}  
\label{th7_2}
\end{equation}
Let us first consider the case in which $\rho_{i-1}(\theta+\delta)>
\rho_{i-1}(\theta-\delta)$ and $\rho_{i}(\theta+\delta)\leq
\rho_{i}(\theta-\delta)$. When $\delta\leq\delta_{0}$, (\ref{th7_2}) ensures that
(\ref{th7_1}) can be rewritten as
\[ \sum_{i=1}^{m}\int_{\rho_{i-1}(\theta-\delta)}^{\rho_{i-1}(\theta+\delta)}
   E_{\xi_{2}}[(L(X(\theta+\delta,\xi))-L(X(\theta-\delta,\xi)))^{2}]d\xi_{1} \]
\[ +\sum_{i=1}^{m}\int_{\rho_{i-1}(\theta+\delta)}^{\rho_{i}(\theta+\delta)}
   E_{\xi_{2}}[(L(X(\theta+\delta,\xi))-L(X(\theta-\delta,\xi)))^{2}]d\xi_{1} \]
\[ +\sum_{i=1}^{m}\int_{\rho_{i}(\theta+\delta)}^{\rho_{i}(\theta-\delta)}
   E_{\xi_{2}}[(L(X(\theta+\delta,\xi))-L(X(\theta-\delta,\xi)))^{2}]d\xi_{1} \]
\[ =\sum_{i=1}^{m}\int_{\rho_{i-1}(\theta-\delta)}^{\rho_{i-1}(\theta+\delta)}
   E_{\xi_{2}}[(L(F_{i-1}^{-1}(\xi_{2}))-L(F_{i}^{-1}(\xi_{2})))^{2}]d\xi_{1} \]
\[\;\;\;+\sum_{i=1}^{m}\int_{\rho_{i}(\theta+\delta)}^{\rho_{i}(\theta-\delta)}
   E_{\xi_{2}}[(L(F_{i+1}^{-1}(\xi_{2}))-L(F_{i}^{-1}(\xi_{2})))^{2}]d\xi_{1} \]
\[ =\sum_{i=1}^{m} E[(L(F_{i-1}^{-1}(\xi)-L(F_{i}^{-1}(\xi)))^{2}]
      \rho_{i-1}'(\theta)2\delta  \]
\[ \;\;\; +\sum_{i=1}^{m}E[(L(F_{i+1}^{-1}(\xi))-L(F_{i}^{-1}(\xi)))^{2}]
       \rho_{i}'(\theta)2\delta +o(\delta)                      \]
By considering every case of $\rho_{i-1}(\theta+\delta)\leq
\rho_{i-1}(\theta-\delta)$ and $\rho_{i}(\theta+\delta)\leq
\rho_{i}(\theta-\delta)$, $\rho_{i-1}(\theta+\delta)>
\rho_{i-1}(\theta-\delta)$ and $\rho_{i}(\theta+\delta)>
\rho_{i}(\theta-\delta)$, and $\rho_{i-1}(\theta+\delta)\leq
\rho_{i-1}(\theta-\delta)$ and $\rho_{i}(\theta+\delta)>
\rho_{i}(\theta-\delta)$, we obtain that
\[    E[(L(X(\theta+\delta,\xi))-L(X(\theta-\delta,\xi)))^{2}]
           = M_{3}(\theta)\delta+o(\delta).  \]
It follows from (\ref{h_var}) that $Var[h] = (1/4)M_{3}(\theta)/\delta+o(1/\delta)$.
Applying Corollary 2, it is easy to see that the convergence rate for
the KW algorithm is $n^{-2/5}$. We have thus completed the proof.
\hspace{0.1in}\rule{1.5mm}{3.0mm}
\vspace{3 mm}

Similarly, we can prove that the rate $n^{-2/5}$ remains valid 
for the case in which the random variable $X(\theta,\xi)=
F_{i}^{-1}(\theta,\xi_{2})$ is generated by setting $\xi_{2}=
(\xi_{1}-\rho_{i-1}(\theta))/p_{i}(\theta)$. 
\vspace{3 mm}

{\sc Theorem 8}. {\em Suppose that $p_{i}(\theta)>0$ is differentiable and the
following function exists and is finite for all $\theta$:
\[      M_{4}(\theta) = \sum_{i=1}^{m}[L(F_{i+1}^{-1}(1^{-}))
                -L(F_{i}^{-1}(0^{+}))]^{2}|\rho_{i}'(\theta)|,         \]
where $F_{i+1}^{-1}(1^{-})=\lim_{\xi\uparrow 1}F_{i+1}^{-1}(\xi)$ and
$F_{i}^{-1}(0^{+})=\lim_{\xi\downarrow 0}F_{i}^{-1}(\xi)$. If 
$M_{4}(\theta)>0$ is bounded for all $\theta$ and (A1)-(A4) are 
satisfied, then the 
convergence rate for the KW algorithm is $n^{-2/5}$}.
\vspace{3 mm}

The previous Theorems 7 and 8 show that the convergence rate for the KW 
algorithm is $n^{-2/5}$ when the composition method is used. This rate
does not depend on how many random numbers are used in the generation 
of random variables. We would emphasize that it is unlikely for each of 
$M_{3}(\theta)$, $M_{4}(\theta)$ to be zero in practice.
\vspace{5mm}

{\bf 5. The MD algorithm with CRN.}\hspace{3mm}
In this section, we examine the rates of convergence for the MD
 algorithm under CRN. As shown in the previous section, the use of
 CRN largely affects $E[h_n^2]$ and thus the reduction of the variance $Var[h_n]$.
 The analysis in the previous Section 4 provides direct information on $E[h_n^2]$. Therefore,
 in this section we directly work $E[h_n^2]$ without going through $Var[h_n]$. We may
 represent $E[h_n^2]$ in the following form.
 \begin{equation}
E[h_n^2] \leq \tilde{c}\delta_n^\gamma.
\label{h3}
\end{equation}
Recall that  $\gamma=-2$ for independent samplings of $X(\theta,\xi)$ without CRN.
With CRN, $\gamma=-1$ if $M_1(\theta)>0$ and  $\gamma=0$ if $M_1(\theta)=0$.
By following the same arguments as in the proof of Theorem 3 and applying (\ref{h3}) directly for  
$E[h_i^2]$ in (\ref{eqn2}), we obtain the following theorem.
\vspace{3 mm}

{\sc Theorem 9.} {\em Assume (B1)-(B3) and (\ref{h3}). Then we have
\begin{equation}
    E[J(\hat{\theta}_{n})-J(\theta^{*})] \leq  \frac{C_1}{n a_n}+\frac{\tilde{C}_2}{n}\sum_{i=1}^n a_i \delta_i^\gamma
+ \frac{C_5}{n}\sum_{i=1}^n \delta_i^\beta,
\label{eqn_master2}
\end{equation}
where $\tilde{C}_2 = \tilde{c}/(2\kappa)$, $C_1$ and $C_5$ are specified in Theorem 3.}

The following Corollary 7 summarizes the rates of convergence for the MD algorithm with 
using CRN in calculating the finite difference (\ref{h_generic}) and (\ref{th5_1}).          
\vspace{3 mm}

{\sc Corollary 7.} {\em Assume (B1)-(B3) and (\ref{h3}). Denote
\[ \tilde{H}_n = \frac{C_1}{n a_n}+\frac{\tilde{C}_2}{n}\sum_{i=1}^n a_i \delta_i^\gamma
+ \frac{C_5}{n}\sum_{i=1}^n \delta_i^\beta.
\]
Then
\begin{itemize}
\item[(i)] if $\gamma=-1$, the best possible rate of convergence for the upper bound $H(n)$ is $n^{-1/3}$
when the one-sided finite difference (\ref{th5_1}) is used,
\item[(ii)]  if $\gamma=-1$, the best possible rate of convergence for the upper bound $H(n)$ is $n^{-2/5}$
when the symmetric finite difference (\ref{h_generic}) is used.
\item[(iii)] if $\gamma=0$, the best possible rate of convergence for the upper bound $H(n)$ is $n^{-1/2}$
when either the one-sided finite difference (\ref{th5_1}) or the symmetric finite difference (\ref{h_generic}) is used,
\end{itemize}
}

\vspace{5mm}

{\bf 6. Generalization and applications}. \hspace{3mm}
In Sections 4-5, all the results are obtained for one dimensional random 
variables only. In this section, we extend the results to a case of 
multivariates, which is not difficult but very tedious. Assume that
 $J(\theta) = E_{\xi}[L(X(\theta,\xi))]$, 
where the multidimensional random variable $X(\theta,\xi)=
[X_{1}(\theta,\xi),X_{2}(\theta,\xi),...,X_{m}(\theta,\xi)]^{T} \in 
R^{m}$. For
each $i$, $X_{i}(\theta,\xi)=X_{i}(\theta,\xi_{i})\in R$, $\xi_{i}$ 
is uniform on $[0, 1)$. We only consider the case in which each 
$X_{i}(\theta,\xi_{i})$ is generated from $\xi_{i}$ using the 
inversion method. To avoid repetition, we list the result without
proof which is very similar to that of Theorem 5.

Assume that $J(\theta)\in R$ and $\theta\in \Theta$. For each $i$, $1\leq i 
\leq m<\infty$, let $F_{i}(\theta,x)$ be the distribution function of 
$X_{i}(\theta,\xi_{i})$ with the decomposition that
\[  \frac{dF_{i}(\theta,x)}{dx}=\left\{ \begin{array}{ll}
      0, & \mbox{if }x\in B_{i,j}^{0}(\theta)=
                [b_{i,j}(\theta),c_{i,j}(\theta)] \\
      f_{i,j}(\theta,x), & \mbox{if }x\in B_{i,j}^{+}(\theta)=
                (c_{i,j}(\theta),b_{i,j+1}(\theta)),
      \end{array}\right.    \]
where $\bigcup_{j}\{B_{i,j}^{0}(\theta)\bigcup B_{i,j}^{+}(\theta)\} = R$
for all $i$, $f_{i,j}(\theta,x) > 0$ for any $x\in B_{i,j}^{+}(\theta)$.
It is possible that $F_{i}(\theta,x)$ is discontinuous at $b_{i,j}(\theta)$.
\vspace{3 mm}

{\sc Theorem 10}. {\em Assume Assumptions (A1)-(A4) and, in addition,
\begin{itemize}
\item[(C1)'.] $L(X)$ is continuously differentiable in $X$, $L(X)$ and 
             $L'_{X_{i}}(X)$ are bounded for all $i$;
\item[(C2)'.] for each $i$,
      \[  \sum_{j}E[(\max_{\theta}\left(\frac{\partial F_{i,j}(\theta,x)}
          {\partial\theta}\right)^{2}/\frac{\partial F_{i,j}(\theta,x)}
           {\partial x})I_{B_{i,j}^{+}(\theta)}]<\infty;  \]
\item[(C3)'.] for all $i, j$, $b_{i,j}(\theta)$ is continuously 
        differentiable in $\theta$, and
      \( \sum_{j}\max_{\theta}(b_{i,j}'(\theta))^{2}< \infty; \)
\item [(C4)'.] for all $i, j$, the functions $F_{i}(\theta,c_{i,j}(\theta))$
             and $F_{i}(\theta,b_{i,j}^{-}(\theta))$ are continuously 
             differentiable in $\theta$, and
      \( \sum_{j}\max_{\theta}|F'_{i}(\theta,c_{i,j}(\theta))|<\infty,\)
   \( \sum_{j}\max_{\theta}|F'_{i}(\theta,b_{i,j}^{-}(\theta))|<\infty,\)
\end{itemize}
Define
\(  \tilde{M}_{1}(\theta) = 
         \sum_{i,j}(L(c_{i,j}(\theta))-L(b_{i,j}(\theta)))^{2}
           |F'_{i}(\theta,c_{i,j}(\theta))|.       \)
Then $\tilde{M}_{1}(\theta)\geq 0$ is bounded. If $\tilde{M}_{1}(\theta)>0$ 
for all $\theta$, the best possible convergence rate for the KW algorithm (2) 
with $h_{n}$ defined by (24) is $n^{-2/5}$. This rate is attained 
by choosing $a_{n}=an^{-1}$, $a>2/(5K_{1})$, and $\delta_{n}=n^{-1/5}$}. 

\vspace{3mm}
Similar results can be obtained if other methods are used in the 
generation of random variables or if $L(X)$ is a piecewise continuous
function of $X$. The analysis can be applied to 
general problems such as Monte Carlo optimization of queueing 
systems and other general systems.
Although such a generalization is not trivial, the basic
idea is the same except that the analysis becomes
tedious and lengthy. Next we illustrate an application of Theorem 10
to the optimization of queueing systems [see, e.g. Kleinrock (1976)].
\vspace{4mm}

{\sc Example 1. GI/G/1 queue with single class of customers}. In a GI/G/1
queue, there is one server (such as a teller in a bank) and one queue. Upon
its arrival, a customer enters the server for service if the server is free,
otherwise it joins the queue and waits for its turn. The service 
discipline is first-come-first-serve. The server cannot be free if there
is at least one customer waiting in the queue. Assume that the 
distribution of interarrival times is $G_{a}(t)$ and the distribution 
of service times is $G_{s}(\theta,t)= p(\theta)G_{s}^{1}(t)+
(1-p(\theta))G_{s}^{2}(t)$. For simplicity, we assume that $G_{a}(t)$, 
$G_{s}^{1}(t)$, and $G_{s}^{2}(t)$ are independent of $\theta$ 
and $\int t^{2}dG_{s}^{j}(t)<+\infty, j = 1, 2$, $p(\theta)$ is 
continuously differentiable in $\theta$, $G_{a}(t), G_{s}^{j}(t),
j=1,2,$ are strictly increasing and continuously differentiable in $t$. 
In queueing theory, 
the system time of a customer is defined as the time period from 
its arrival till departure. Let $L(X(\theta,\xi))$ be the average 
system time of the first N customers
\[  L(X(\theta,\xi)) = \frac{1}{N}\sum_{i=1}^{N}T_{i}(\theta,\xi),  \]
where $T_{i}(\theta,\xi)$ is the system time of the $i$th customer.
Then $J(\theta)=E[L(X(\theta,\xi))]$ is the mean system time of the 
first $N$ customers. We want to find the optimal parameter $\theta^{*}$ 
to minimize $J(\theta)$. It is known that the
analytical form of $J(\theta)$ is not available for general $G_{a}(t)$,
$G_{s}^{1}(t)$, and $G_{s}^{2}(t)$ [e.g. Kleinrock (1976)]. So we find 
$\theta^{*}$ via the KW algorithm. Assume that the queue is initially
empty. According to Lindley's equation [e.g. Kleinrock (1976)]:
\begin{equation}
    T_{i}(\theta,\xi) =
      \max\{T_{i-1}(\theta,\xi)-A_{i},0\}
                 +S_{i},\;\;   T_{0}(\theta,\xi) = 0,  
\label{queue1}
\end{equation}
where $A_{i}$ is the interarrival time between the $(i-1)$th and 
the $i$th customer, $S_{i}$ is the service time
of the $i$th customer. The distributions of $A_{i}$ and 
$S_{i}$ are respectively $G_{a}(t)$ and $G_{s}(\theta,t)$.
We consider two scenarios.

{\em Case 1}. We find $\theta^{*}$ through computer simulation. We
write a program to simulate the GI/G/1 queue. At the $n$th iteration, 
we perform two experiments with the same $\xi_{n}$ to obtain a $h_{n}$
that is defined by (\ref{h_generic}). Consider that the inversion method is used 
in the generation of random variables $A_{i}=G_{a}^{-1}(u_{i}), 
S_{i}=G_{s}^{-1}(\theta,v_{i}), i=1,2,..., N$. Define the random
factor as $\xi=
[u_{1},u_{2},...,u_{N},v_{1},v_{2},...,v_{N}]^{T}$, $A(\xi)= 
[A_{1},A_{2},...,A_{N}]$, $S(\theta,\xi)=[S_{1},S_{2},...,S_{N}]$, 
and $X(\theta,\xi)=[A(\xi), S(\theta,\xi)]^{T}$. Since the function 
$\max\{x,0\}$ is continuous in $x$, $L(X)$ is continuous in $X$.
According to Theorem 10, we know that $\tilde{M}_{1}(\theta)=0$ 
since both $G_{a}(t)$ and $G_{s}(\theta,t)$ are strictly increasing 
and continuously differentiable in $t$. 
Note that $L(X)$ is not differentiable in $X$. However, $L(X)$ is
left and right differentiable with bounded one-sided derivatives.
A simple modification of the proof of Corollary 3 shows that 
$Var[h_{n}]=O(1)$. Therefore, the convergence
rate for the KW algorithm is $n^{-1/2}$. If the composition method is
used in the generation of $S(\theta,\xi)$ according to the distribution
$G_{s}(\theta,t)$, then from Theorems 7 and 8
(which is applicale to the case of multivariates) we know 
that the rate of convergence is $n^{-2/5}$.

{\em Case 2}. Assume that this is a real system and we want to perform
on-line parameter adjustment. Let visualize $n$ as the $n$th day of
service. Suppose that the server serves more than $N$ customers
 each day. At the $n$th day, the server serves customers with 
parameter value $\theta_{n}$ and simultaneously collects information 
of $X(\theta_{n},\xi_{n})$ which simply is a record of interarrival 
times $\{A_{i}^{n}\}$ and service times $\{S_{i}^{n}\}$. 
At the end of the $n$th day, the server calculates
\[   v_{i}=G_{s}(\theta_{n},S_{i}^{n}), i=1,2,3, ..., N.\]
It is easy to verify that each $v_{i}$ is uniform on $[0, 1)$.
Then the server defines $\xi_{n}$ from the preceding $v_{i},
i=1,2, ..., N$, takes a $\delta_{n}>0$, and 
\[ S(\theta_{n}+\delta_{n},\xi_{n})=[S_{1}^{n,1},S_{2}^{n,1},...,S_{N}^{n,1}],
      S_{i}^{n,1}
          =G_{s}^{-1}(\theta_{n}+\delta_{n},G_{s}(\theta_{n},S_{i}^{n})), 
       \;\;  i=1,2,..., N;   \]
\[ S(\theta_{n}-\delta_{n},\xi_{n})=[S_{1}^{n,2},S_{2}^{n,2},...,S_{N}^{n,2}],
      S_{i}^{n,2}
          =G_{s}^{-1}(\theta_{n}-\delta_{n},G_{s}(\theta_{n},S_{i}^{n})),
        \;\; i=1,2, ..., N. \]
If $G_{s}(\theta,t)=1-e^{-t/\theta}$ is exponential, then $S_{i}^{n,1}=
(\theta_{n}+\delta_{n})S_{i}^{n}/\theta_{n}$, $S_{i}^{n,2}=
(\theta_{n}-\delta_{n})S_{i}^{n}/\theta_{n}$. 
With the values of $A(\xi_{n}), S(\theta_{n}+\delta_{n},\xi_{n}), 
S(\theta_{n}+\delta_{n},\xi_{n})$, from (40) and the form of 
$L(X(\theta,\xi))$, the server computes $L(X(\theta_{n}+
\delta_{n},\xi_{n}))$ and $L(X(\theta_{n}-\delta_{n},\xi_{n}))$, which
determines a $h_{n}$. With this $h_{n}$, the server updates the parameter
$\theta_{n+1}$ according to the KW algorithm (\ref{kw}) for the next $(n+1)$th
day. In
such a way, we have formulated an on-line optimization problem that mimics
the Monte Carlo optimization. Its convergence can be analyzed similarly
to that of Case 1. Our purpose here is simply to point out that the results
of this paper are not restricted to Monte Carlo optimization.

\vspace{5mm}

{\bf 7. Summary}. \hspace{3mm}
So far, we have examined several variations of the KW algorithm and 
the MD algorithm under
the symmetric finite difference, the one-sided finite difference, and
the use of CRN when different methods are used in the generation of 
random variables. The results of 
this paper, together with previous results on the KW algorithm without 
the use of CRN [c.f. Fabian (1971); Kushner and Clark (1978)], provide a 
complete view toward the rates of convergence for the KW algorithm.
For the ease of comparison, we summarize all the results in the 
following table. 
\vspace{5 mm}
%\newpage

\centerline{Table I. Rates of convergence for the KW/MD algorithm}
\vspace{5 mm}

\begin{tabular}{|c|c|c|c|c|}  \hline
    & \(\begin{array}{c} \mbox{with CRN}\\ \mbox{$h$(\ref{h_generic})}
      \end{array}\) & \(\begin{array}{c} \mbox{with CRN}\\ 
      \mbox{$h$(\ref{th5_1})} \end{array}\) & \(\begin{array}{c} 
      \mbox{without CRN}\\ \mbox{$h$(\ref{kw_h})}
      \end{array}\) & \(\begin{array}{c} \mbox{without CRN}\\ 
      \mbox{$h$(\ref{kw_1sideh})} \end{array}\) \\ \hline
\(\begin{array}{c} \mbox{inversion:} \\ M_{1}(\theta)\neq 0 \end{array}\) & 
    $n^{-2/5}$ & $n^{-1/3}$ & $n^{-1/3}$ & $n^{-1/4}$  \\ \hline
\(\begin{array}{c} \mbox{inversion:} \\ M_{1}(\theta)= 0 \end{array}\) & 
    $n^{-1/2}$ & $n^{-1/2}$ & $n^{-1/3}$ & $n^{-1/4}$ \\ \hline
\(\begin{array}{c} \mbox{rejection:} \\ \mbox{general} \end{array}\) 
& $n^{-2/5}$ & $n^{-1/3}$ & $n^{-1/3}$ & $n^{-1/4}$  \\ \hline
\(\begin{array}{c} \mbox{composition:} \\ \mbox{general} \end{array}\) 
 & $n^{-2/5}$ & $n^{-1/3}$ & $n^{-1/3}$ & $n^{-1/4}$  \\ \hline
\end{tabular}
\vspace{10 mm}

In Table I, $h$(\ref{kw_h}), $h$(\ref{kw_1sideh}), $h$(\ref{h_generic}), and $h$(\ref{th5_1}) refer to the 
finite-difference approximation $h_{n}$ defined by (\ref{kw_h}), (\ref{kw_1sideh}), (\ref{h_generic}), 
and (\ref{th5_1}), respectively. The phrase ``without CRN'' refers to using 
independent samples in calculating the finite difference $\{h_n\}$,
which excludes sampling schemes that may lead to correlations between 
the samples. In other words, ``without CRN'' simply means that $\xi_{1,n}$ and $\xi_{2,n}$ are 
independent in (\ref{kw_h}) and (\ref{kw_1sideh}). When the inversion method is used in 
the generation of random variables
and when $M_{1}(\theta)=0$, we assume that Corollaries 3 and 4
are applicable. Results pertaining to the KW algorithm without the use 
of CRN can be found in, for example, Fabian (1971), 
and Kushner and Clark (1978).
\vspace{3 mm}

Generally speaking, the
use of CRN is always helpful in accelerating the convergence of the 
KW algorithm or the MD algorithm. In some cases, such as when $M_{1}(\theta)=0$ in Theorem 5,
CRN helps a lot. In some of other cases, CRN may help less much. 
When the inversion method is used and when $M_{1}(\theta)=0$,
the convergence rate can reach the best possible rate for the two types 
of stochastic approximation algorithms. The remark at the end of 
Subsection 3.2 shows that, as far as the convergence rate of the KW 
algorithm is concerned, the inversion method is superior to 
the rejection method. Note that inversion can also be used to generate
random variables with distributions of the form 
$\sum_{i}p_{i}(\theta)F_{i}(\theta,x)$. A comparison of Theorem 4 and
Theorems 7 and 8 shows that inversion is also superior to composition.
When the distribution function $F(\theta,x)$ of $X(\theta,\xi)$
is strictly increasing and continuous, a close examination of the 
inversion, rejection, and composition methods shows that
$X(\theta,\xi)$ is continuous in $\theta$ if it is generated from
inversion. However, $X(\theta,\xi)$ is discontinuous in $\theta$ if
it is generated from either rejection or composition. It is such 
a distinction of continuity that determines the rates of
the convergence for the KW algorithm. 
%As a general guidance, one 
%should use (continuous) inversion whenever it is possible.

\vspace{5mm}

\centerline{REFERENCES}
\begin{enumerate}
\item {\sc P. Bratley, B. Fox, and L. Schrage,} {\em A Guide to 
      Simulation}, Springer-Verlag, New York, 1983
\item {\sc D.L. Burkholder,} {\em On a class of stochastic approximation
      processes,} Annals of Mathematical Statistics, 27 (1956), pp. 1044-1059.
\item {\sc S. Cambanis, G. Simons, and W. Stout,} {\em Inequalities for 
      $Ek(X,Y)$ when marginals are fixed,} Z. Whar. Geb. 36 (1976),
      pp. 285-294.
\item {\sc K.L. Chung,} {\em On a stochastic approximation method,} 
      Annals of Mathematical Statistics, 25 (1954), pp. 463-483.
\item {\sc R. W. Conway,} {\em Some tactical problems in digital simulation,} 
      Management Science, 10 (1963), pp. 47-61.
\item {\sc L. Devroye,} {\em Coupled samples in simulation,} Operations 
      Research, 38 (1990), pp. 115-126.
\item {\sc V. Dupa$\check{c}$,} {\em On the Kiefer-Wolfowitz approximation 
      method,} Casopis Pest. Math, 82 (1957), pp. 47-75.
 \item {\sc J.C. Duchi, A. Agarwal, M. Johansson, and M.I. Jordan,}
      {\em  Ergodic Mirror Descent,} SIAM Journal on Optimization, 22 (2012), pp. 1549-1578. 
 \item {\sc J. Duchi, M.I. Jordan, M. Wainwright, and A. Wibisono,} {\em 
      Finite sample convergence rates of zero-order stochastic optimization methods,} 
      In {\em Advances in Neural Information Processing Systems (NIPS)}, 
      P. Bartlett, F. Pereira, L. Bottou and C. Burges (Eds.), 2013. 
 %\item {\sc J. Duchi, M. I. Jordan, M. Wainwright, and A. Wibisono,} {\em
 %     Optimal rates for zero-order optimization: the power of two function 
 %     evaluations,} arXiv:1312.2139, 2013b.       
\item {\sc A. Dvoretzky,} {\em On stochastic approximation,} Proc. Third
      Berkeley Symp. Math. Statis. Prob., 1 (1956), pp. 39-56.
\item {\sc V. Fabian,} {\em Stochastic approximation,} In {\em Optimizing
      Methods in Statistics}, J.S. Rustagi (ed.), Academic Press, New York, 1971.
\item {\sc G.S. Fishman,} {\em Correlated simulation experiments,}
      Simulation, 23 (1974), pp. 177-180.
\item {\sc W.R. Franta,} {\em The Process View of Simulation}. North
      Holland, New York, 1975.
\item {\sc P. Glasserman and D. Yao,} {\em  Some guidelines and guarantees 
      for common random numbers,} Management Sciences, 38(1992), pp. 884-908.
%\item {\sc Glynn, P.W.} 1985. Regenerative structure of Markov chains simulated
%      via common random numbers. {\em O.R. Letters} 4 : 49-53.
\item {\sc J.M. Hammersley and D.C. Handscomb,} {\em Monte Carlo
      Methods}, Methuen, London, 1964.
%\item {\sc Heidelberger, P., and D.L. Iglehart.} 1979. Comparing stochastic 
%      systems using regenerative simulation and common random 
%      numbers. {\em Advances in Applied Probability} 11 : 804-819.
\item {\sc R.G. Heikes, D.C. Montogomery, and R.L. Rardin,} {\em Using
      common random numbers in simulation experiments,} Simulation,
      27 (1976), pp. 81-85.
\item {\sc Y.C. Ho and X.R. Cao,} {\em Perturbation Analysis of
      Discrete Event Dynamic Systems,} Kluwer Academic Publishers, Boston, 1991.
\item {\sc H. Kesten,} {\em Accelerated stochastic approximation,} 
      Annals of Mathematical Statistics, 29 (1958), pp. 41-59. 
\item {\sc J. Kiefer and J. Wolfowitz.} {\em Stochastic estimation of the
      maximum of a regression function,} Annals of Mathematical 
      Statistics, 23 (1952), pp. 462-466.
\item {\sc J.P.C. Kleijnen,} {\em Statistical Techniques in Simulation,}
      Marcel Dekker, New York, 1974.
\item {\sc L. Kleinrock,} {\em Queueing Systems}, Vol.I, Wiley, New York, 1976.
\item {\sc S.G. Krantz,} {\em Real Analysis and Foundations}, CRC Press,
      Bocan Raton, FL, 1991.
\item {\sc H.J. Kushner and D.S. Clark,} {\em Stochastic Approximation
      Methods for Constrained and Unconstrained Systems}, Springer-Verlag.
      New York, 1978.
\item {\sc A.M. Law and W.D. Kelton.}{\em Simulation Modeling 
      and Analysis}, McGraw-Hill, New York, 1982.
%\item {\sc Mitchell, B.} 1973. Variance reduction by antithetic variates in 
%      GI/G/1 queueing simulation. {\em Operations Research} 21 : 988-997.
\item {\sc A. Nemirovski, A. Juditsky, G. Lan, and A. Shapiro,} {\em 
      Robust stochastic approximation approach to
      stochastic programming,} SIAM Journal on Optimization, 19 (2009), pp. 1574-1609.
\item {\sc A. Nemirovski, and D. Yudin,} {\em Problem Complexity and 
         Method Efficiency in Optimization}, Wiley, New York, 1983.
%\item {\sc Reiman, M.I., and A. Weiss.} 1986. Sensitivity analysis via 
%      likelihood ratios. {\em Proc. of the 1986 Winter Simulation 
%      Conference}, 285-289.
\item {\sc H. Robbins, and S. Monro,} {\em A stochastic approximation method,}
      Annals of Mathematical Statistics, 22 (1951), pp. 400-407.
%\item {\sc Rubinstein, R.Y., and G. Samorodnitsky.} 1985. Variance reduction
%      by the use of common and antithetic random variables. 
%      {\em J. of Statistics and Computer Simulation} 22 : 161-180.
\item {\sc J. Sacks,} {\em Asymptotic distribution of stochastic approximation
      procedures,} Annals of Mathematical Statistics, 29 (1958), pp. 373-405.
\item {\sc A.N. Shiryayev,} {\em Probability}, Springer-Verlag, 
      New York, 1984.
\item {\sc M.T. Wasan,} {\em Stochastic Approximation,} Cambridge
      University Press, Cambridge, England, 1969.
\item {\sc W. Whitt,} {\em Bivariate distributions with given marginals,}
      Ann. Math. Stat., 4 (1976), pp. 1280-1289.
%\item {\sc Wright, R.D., and T.E. Ramsay, Jr.} 1979. On the effectiveness of
%      common random numbers. {\em Management Sciences} 25 : 649-656.
\end{enumerate}
\vspace{10 mm}

%\hspace{3.0in}\parbox{3.4in}{{\sc Department of Systems Science 
%
%\hspace{0.5in} and Mathematics 
%
%Washington University 
%
%St. Louis, MO 63130}}
\end{document}